\documentclass[11pt,reqno]{amsart}

\usepackage{fullpage}
\usepackage{amssymb,amsmath,amsthm,latexsym}
\usepackage[usenames]{color}
\usepackage{graphicx}
\definecolor{webgreen}{rgb}{0,.5,0}
\definecolor{webbrown}{rgb}{.6,0,0}
\usepackage{tabularx}
\usepackage{array}
\usepackage{mathtools}
\usepackage{mathrsfs}
\usepackage{enumerate}
\usepackage{float}
\usepackage[colorlinks=true,
linkcolor=webgreen,
filecolor=webbrown,
citecolor=webgreen]{hyperref}

\DeclarePairedDelimiter\floor{\lfloor}{\rfloor}
\newtheorem{theorem}{Theorem}[section]
\newtheorem{conjecture}[theorem]{Conjecture}
\newtheorem{corollary}[theorem]{Corollary}
\newtheorem{definition}[theorem]{Definition}

\newtheorem{lemma}[theorem]{Lemma}
\newtheorem{question}[theorem]{Question}
\newtheorem{proposition}[theorem]{Proposition}
\theoremstyle{remark}
\newtheorem{remark}[theorem]{Remark}

\newcommand{\cA}{\mathcal{A}}
\newcommand{\cB}{\mathcal{B}}
\newcommand{\tr}{\operatorname{tr}}

\begin{document}
\title{Integral Means of Fibonacci Numbers: Divisibility, Density, and Arithmetic Structure}
\author[Sagar Mandal]{Sagar Mandal$^*$}
\address{Department of Mathematics\\Indian Institute of Technology Ropar\\Punjab, India.}
\email{sagar.25maz0008@iitrpr.ac.in, sagarmandal31415@gmail.com}
\thanks{$^*$The first author is financially supported by the University Grants Commission (UGC), Government of India, through the award of the Junior Research Fellowship (JRF) with ref. No. 241620111598.}
\author[Mohan]{Mohan$^{\dagger}$}
\address{ School of Engineering and Technology, Vivekananda Institute of Professional Studies - Technical Campus, AU- Block (Outer Ring Road), Pitampura, Delhi - 110034, India}
\email{mohan98math@gmail.com}
\thanks{$^{\dagger}$The corresponding author.}

\maketitle
\begin{flushright}
\textit{\textbf{-- In memory of Prof. Florian Luca}}\\[2mm]
\textit{\textbf{(1969--2026)}}
\end{flushright}
\begin{abstract}
Motivated by recent works of Fatehizadeh, Yaqubi, Flynn-Connolly, Luca and Marques on the arithmetic mean of Fibonacci numbers, we study divisibility properties of sums of Fibonacci numbers whose indices form arithmetic progressions. We show that infinitely many even and odd $n$ divide such sums, and also the corresponding sums of squares and cubes, and we find all $n$ for which the geometric and harmonic means of $F_1,\ldots,F_n$ are integers. Further, we study the set $\cA$ of integers $n\geq2$ for which the arithmetic mean of $F_1,\ldots,F_n$ is an integer. We show that $\cA$ is both rich and sparse: it contains arbitrarily long arithmetic progressions and does not contain odd prime powers, yet $\#\{n\leq x: n\in\cA\}\ll x/\sqrt{\log x}$. This density-zero bound holds more generally for every set $\{n: n\mid\sum_{i=0}^{n-1}F_{r+id}\}$. We also establish necessary conditions for three consecutive integral means (a question of Luca and Marques). 
\end{abstract}

\section{Introduction}
The Fibonacci sequence $\{F_n:n\geq0\}$ (\href{https://oeis.org/A000045}{A000045}) is the binary recurrence sequence defined as: $F_0=0,F_1=1,$ and $ F_{n+2}=F_{n+1}+F_n$ for $n\geq0$. The companion Lucas sequence $\{L_n:n\geq 0\}$ (\href{https://oeis.org/A000032}{A000032}) is defined by $L_0=2$, $L_1=1$ and $L_{n+2}=L_{n+1}+L_n$. These are among the most extensively studied sequences in number theory, see \cite{new1,r1,r2,r10,new2,Kom20,Kom25,Kom26,r3,new3,new4,new5}. One classical problem is to determine when certain expressions involving Fibonacci numbers are integers. In recent years, several authors have investigated the arithmetic mean of Fibonacci numbers (see \cite{r1,r2,r3,r12}). For example, Fatehizadeh and Yaqubi \cite{r1} studied the arithmetic mean of the first $n$ Fibonacci numbers, which is defined by
$$A_F(n):=\frac{1}{n}\sum_{i=1}^nF_i.$$
They proved that there are infinitely many even integers $n$ such that $A_F(n)\in\mathbb{N}$ (see \cite{r12} for generalizations). In the same paper, they gave the following conjectures:
\begin{conjecture}
There exist infinitely many odd positive integers $n$ such that
$A_F(n)\in\mathbb{N}.$
\end{conjecture}

\begin{conjecture}\label{Conjecture1}
There exist infinitely many pairs of consecutive positive integers
$(n,n+1)$ such that both
$A_F(n),~ A_F(n+1)$
are integers.
\end{conjecture}
The first conjecture was later proved by Flynn-Connolly \cite{r2}, whereas the second was established by Luca and Marques \cite{r3}.

Motivated by the arithmetic mean of the first $n$ Fibonacci numbers, we are interested in the arithmetic mean of the first $n$ Fibonacci numbers whose
indices form arithmetic progressions, in particular those with common difference $2$ and $3$.

\begin{theorem}\label{Theorem 2}
\noindent
\begin{enumerate}
    \item[\upshape(1)]
    $\displaystyle{\frac{1}{n}\sum_{i=1}^{n}F_{2i}\in \mathbb{N}}$,
    \item[\upshape(2)] $\displaystyle{\frac{1}{n}\sum_{i=1}^{n}F_{2i-1}\in \mathbb{N}}$,
    \item[\upshape(3)] $\displaystyle{\frac{1}{n}\sum_{i=1}^{n}F_{3i}\in \mathbb{N}}$,
    \item[\upshape(4)] $\displaystyle{\frac{1}{n}\sum_{i=1}^{n}F_{3i-1}\in \mathbb{N}}$,
    \item[\upshape(5)]$\displaystyle{\frac{1}{n}\sum_{i=1}^{n}F_{3i-2}\in \mathbb{N}}$,
\end{enumerate}
for infinitely many odd and even positive integers $n$.
\end{theorem}
We next generalize the above results by considering Fibonacci numbers with indices in arithmetic progression and, subsequently, their alternating sums.

\begin{theorem}\label{thm:lambda}
Let $r\geq0$ and $d\geq1$ be integers, and let $\lambda\in\{-1,1\}$. Then there exist infinitely
many even positive integers $n$ such that
\[
n \mid \sum_{i=1}^{n}\lambda^{i}F_{r+di}.
\]
Moreover, if $\lambda=1$ and $4\mid d$, then $n$ may additionally be chosen odd.
\end{theorem}
\begin{corollary}\label{Theorem 1}
Let $r \geq 0$ and $d \geq 1$ be integers. Then
$$\frac{1}{n} \sum_{i=1}^{n}F_{r+di}\in \mathbb{N}$$
for infinitely many even positive integers $n$.
Further, if $4\mid d$, then there are infinitely many odd $n$ such that
$$\frac{1}{n} \sum_{i=1}^{n}F_{r+di}\in \mathbb{N}.$$
\end{corollary}

\begin{remark}
In general, it is unknown whether there exist infinitely many odd integers $n$ such that $n\mid \displaystyle{\sum_{i=1}^{n}F_{r+di}},$ whenever $4 \nmid d$ with $d\geq5$.     
\end{remark}

The preceding results lead us to consider a more general class of divisibility problems involving powers of Fibonacci numbers.\\
We start with the following theorems, which concern the arithmetic mean of the squares of Fibonacci numbers whose indices form arithmetic progressions with common differences $1$, $2$ and $3$.

\begin{theorem}\label{Theorem 5}
\noindent
\begin{enumerate}
    \item[\upshape(1)] $\displaystyle{\frac{1}{n}\sum_{i=1}^{n}F_{i}^2}\in \mathbb{N}$,
    \item[\upshape(2)] $\displaystyle{\frac{1}{n}\sum_{i=1}^{n}F_{2i}^2}\in \mathbb{N}$,
    \item[\upshape(3)] $\displaystyle{\frac{1}{n}\sum_{i=1}^{n}F_{2i-1}^2}\in \mathbb{N}$,
    \item[\upshape(4)] $\displaystyle{\frac{1}{n}\sum_{i=1}^{n}F_{3i}^2}\in \mathbb{N}$,
\end{enumerate}
for infinitely many even and odd positive integers $n$.
\end{theorem}
\begin{theorem}\label{theorem p}\noindent
\begin{enumerate}
\item[\upshape(1)] $\displaystyle{\frac{1}{n}\sum_{i=1}^{n}F_{3i-1}^2}\in \mathbb{N}$,
\item[\upshape(2)] $\displaystyle{\frac{1}{n}\sum_{i=1}^{n}F_{3i-2}^2}\in \mathbb{N}$,
\end{enumerate}
for infinitely many even positive integers $n$.
Moreover, for every odd prime $p$,
$$\frac{1}{p}\sum_{i=1}^{p} F_{3i-1}^{2}\not\in \mathbb{N}\quad \text{and} \quad \frac{1}{p} \sum_{i=1}^{p} F_{3i-2}^{2}\not\in \mathbb{N}.$$
\end{theorem}
\noindent Therefore, no odd prime has the divisibility property for the two sums in Theorem \ref{theorem p}. Moreover, a computer search shows that there is no odd $n$ with $1<n<6\cdot 10^4$ (see \cite{Codes} for Python codes) dividing either of these two sums. This computational evidence leads us to propose the following conjecture.
\begin{conjecture}\label{conj:odd3}\noindent
\begin{enumerate}
\item[\upshape(1)] $\displaystyle{\frac{1}{n}\sum_{i=1}^{n}F_{3i-1}^2}\in \mathbb{N}$,
\item[\upshape(2)] $\displaystyle{\frac{1}{n}\sum_{i=1}^{n}F_{3i-2}^2}\in \mathbb{N}$,
\end{enumerate}
for only finitely many odd positive integers $n>1$.
\end{conjecture}
In the next result, we present an analogous result to Corollary \ref{Theorem 1}.

\begin{theorem}\label{Theorem 4}
Let $r \geq 0$ and $d \geq 1$ be integers. Then there exist infinitely many even positive integers $n$ such that
$$n\mid \sum_{i=1}^{n}F_{r+di}^2.$$
Further, if $4\mid d$, then there are infinitely many odd $n$ with the same property.
\end{theorem}

The following theorems concern the arithmetic mean of the cubes of Fibonacci numbers whose indices form arithmetic progressions with common differences $1$, $2$ and $3$.

\begin{theorem}\label{Theorem 8}\noindent
\begin{enumerate}
    \item[\upshape(1)] $\displaystyle{\frac{1}{n}\sum_{i=1}^nF_i^3}\in \mathbb{N},$
    \item[\upshape(2)] $\displaystyle{\frac{1}{n}\sum_{i=1}^nF_{2i}^3}\in \mathbb{N}$,
    \item[\upshape(3)] $\displaystyle{\frac{1}{n}\sum_{i=1}^nF_{2i-1}^3}\in \mathbb{N}$,
\end{enumerate}
for infinitely many even and odd positive integers $n$.
 Furthermore, $\displaystyle{\frac{1}{p}\sum_{i=1}^pF_i^3}\not\in\mathbb{N}$ for every odd prime $p$.
\end{theorem}
\begin{theorem}\label{Theorem 9}\noindent
\begin{enumerate}
    \item[\upshape(1)] $\displaystyle{\frac{1}{n}\sum_{i=1}^nF_{3i}^3}\in \mathbb{N},$
    \item[\upshape(2)] $\displaystyle{\frac{1}{n}\sum_{i=1}^nF_{3i-1}^3}\in \mathbb{N}$,
    \item[\upshape(3)] $\displaystyle{\frac{1}{n}\sum_{i=1}^nF_{3i-2}^3}\in \mathbb{N}$,
\end{enumerate}
for infinitely many even positive integers $n$.
\end{theorem}
\begin{remark}
 Odd positive integers $n>1$ with the divisibility property discussed in Theorem \ref{Theorem 9} for the individual sums also exist, for instance, for $n=39,63,69,207,279$, $n\mid \displaystyle{\sum_{i=1}^{n}F_{3i}^3}$, for $n=3,19$ $n\mid\displaystyle{\sum_{i=1}^{n}F_{3i-1}^3}$, and $7\mid \displaystyle{\sum_{i=1}^{7}F_{3i-2}^3}$. We do not know whether there are infinitely many of them.
\end{remark}
In the following, we establish an analogous result to Corollary \ref{Theorem 1} and Theorem \ref{Theorem 4}.

\begin{theorem}\label{Theorem 7}
Let $r \geq 0$ and $d \geq 1$ be integers. Then there exist infinitely many even positive integers $n$ such that
$$n\mid \sum_{i=1}^{n}F_{r+di}^3.$$
Further, if $4\mid d$, then there are infinitely many odd $n$ with the same property.
\end{theorem}

The divisibility results established in \cite{r1,r2,r3} naturally lead to a more general question. We consider the following question, which for $k=1$ reduces to the problem studied therein.

\begin{question}\label{Q:k} For given $k\geq 1$, are there infinitely many even positive integers and infinitely many odd positive integers $n$ such that  $\displaystyle{\sum_{i=1}^{n}i^{k-1}}$ divides $\displaystyle{\sum_{i=1}^nF_{i}^{k}}$?
\end{question}
\noindent In the following theorem, we answer the preceding question for $k=2$.
\begin{theorem}\label{Theorem 4a}
There exist infinitely many even positive integers and infinitely many odd positive integers $n$ such that
$\displaystyle{\sum_{i=1}^{n}i}$ divides $\displaystyle{\sum_{i=1}^nF_{i}^{2}}$.
\end{theorem}

Motivated by the preceding arithmetic mean problem, we next investigate analogous questions concerning the geometric and harmonic means of the first $n$ Fibonacci numbers. Our main theorems in this direction are as follows.
\begin{theorem}\label{Theorem 10} For $k\geq 1$,
    $\displaystyle{(F_1^kF_2^k \cdots F_n^k)^{\frac{1}{n^k}} \in \mathbb{N}}$ holds only for $n=1,2$.
\end{theorem}

\begin{theorem}\label{Theorem 11} 
\begin{enumerate}
    \item[\upshape(1)] $\dfrac{n}{\displaystyle{\frac{1}{F_1}+\frac{1}{F_2}+\cdots+\frac{1}{F_n}}} \in \mathbb{N}$
     holds only for $n=1,2$.
     \item[\upshape(2)]  $\dfrac{n^2}{\displaystyle{\frac{1}{F_1^2}+\frac{1}{F_2^2}+\cdots+\frac{1}{F_n^2}}} \in \mathbb{N}$
     holds only for $n=1,2,3$.
\end{enumerate}
\end{theorem}
\noindent The $2$-adic method of the proof of Theorem \ref{Theorem 11} does not extend to higher powers (see the discussion before its proof in Section \ref{sec:proofs}). We propose the following conjecture (for $k=3,4,5$ and $n\leq 10^4$, see \cite{Codes} for Python codes).
\begin{conjecture}\label{conj:harm}
Let $k\geq 3$. Then
    $\dfrac{n^k}{\displaystyle{\frac{1}{F_1^{k}}+\frac{1}{F_2^{k}}+\cdots+\frac{1}{F_n^{k}}}} \in \mathbb{N}$ if and only if  $n=1,2$.
\end{conjecture}

Luca and Marques \cite{r3} proved that there are infinitely many $n$ such that $A_F(n)$ and $A_F(n+1)$ are both integers. For three consecutive integers the following question has been asked.

\begin{question}[{\cite[Luca and Marques]{r3}}]\label{conj:Luca}
Does there exist a positive integer $n$ such that $A_F(n)$, $A_F(n+1)$ and $A_F(n+2)$ are all integers?
\end{question}
\noindent More generally, this motivates the following question.

\begin{question}
Let
$$\cA=\Bigl\{n\geq 2:\ n\mid \sum_{i=1}^nF_i\Bigr\}.$$
Does $\cA$ contain arbitrarily long arithmetic progressions? In particular, are there integers $n$ with $n,n+1,n+2\in\cA$?
\end{question}
\noindent Observe that Luca and Marques \cite{r3} essentially proved that there are infinitely many $n$ such that $n,n+1\in\cA$, and there is no  $n$ such that $n,n+1,n+2,n+3\in \cA$. The next result gives a simple sufficient condition for membership in $\cA$ in terms of the Pisano period $\pi(n)$ of $n$, which is the period of the sequence $(F_m\bmod n)_{m\geq0}$ (see Section \ref{sec:prelim}), and shows that $\cA$ is rich in arithmetic structure.
\begin{theorem}\label{thm:Api}Let $n$ be a positive integer. Then the following hold.
\begin{enumerate}
\item[\upshape(1)] If $n\geq2$ and $\pi(n)\mid n$, then $n\in\cA$. In particular, $2^a3^b5^c\in\cA$ for all integers $a\geq3$, $b\geq1$ and $c\geq0$.
\item[\upshape(2)]For every integer $r\geq3$, the numbers $3\cdot2^r, 6\cdot2^r, 9\cdot2^r$ form a $3$-term arithmetic progression contained in $\cA$.
\item[\upshape(3)] For every $k\geq3$, the set $\cA$ contains infinitely many $k$-term arithmetic progressions. More precisely, there exist infinitely many positive integers $N$ such that $jN\in\cA$ for all $j\in\{1,2,\ldots,k\}$.
\end{enumerate}
\end{theorem}
\noindent On the other hand, the elements of $\cA$ satisfy strong congruence restrictions. The next result describes the structure of the even and odd elements of $\cA$, and it will be used to study three consecutive integral means.
\begin{theorem}\label{thm:Aeven}
\noindent
\begin{enumerate}
    \item[\upshape(1)] If $n\in\cA$ is even and $n>2$, then
$
n\equiv0,2,\text{ or }14\pmod{24}.$

\item[\upshape(2)] If $n\in\cA$ is odd, then $5\nmid n$ and $n$ is not a prime power, in particular, $n$ is not a prime. Further, if $3\mid n$, then
    $n\equiv5\text{ or }7\pmod8$,  in particular, $n\not\equiv3,9\pmod{24}$.
\end{enumerate}
\end{theorem}
\noindent As an application of Theorem \ref{thm:Aeven}, in the following theorem, we give a necessary condition such that $A_F(n),A_F(n+1),A_F(n+2)\in\mathbb{N}$.
\begin{theorem}\label{thm:three}
Suppose that $n\geq 1$ and $A_F(n),A_F(n+1),A_F(n+2)\in\mathbb{N}$. Then $n\equiv 0,13$ or $23\pmod{24}$, and there exists an even integer $e\in\{n,n+1\}$, $e\geq4$, such that
$$e\in\cA,\qquad e+1\in\cA,\qquad e\equiv 0\ \text{or}\ 14\pmod{24}.$$
\end{theorem}

We now turn to the Fibonacci numbers that belong to $\cA$.
\begin{definition}[\cite{r2}]
A positive integer $m$ is called \emph{self-summable} if
$$F_m \mid \sum_{i=1}^{F_m} F_i,$$
and in this case $F_m$ is called a \emph{self-summable Fibonacci number}.
\end{definition}
\noindent Since $F_1=F_2=1$, the integers $m=1,2$ are trivially self-summable, and for $m\geq3$, $m$ is self-summable if and only if $F_m\in\cA$. The set $\cA\cup\{1\}$ is the OEIS sequence \href{https://oeis.org/A111035}{A111035}, so the self-summable Fibonacci numbers are exactly the Fibonacci numbers in \href{https://oeis.org/A111035}{A111035}, and the self-summable integers $m$ with $F_m$ odd are listed in \href{https://oeis.org/A381053}{A381053}. Flynn-Connolly \cite{r2} proved the following theorem.
\begin{theorem}\label{thm:FC}
Let $m=2p$ or $m=4p$, where $p$ is an odd prime such that
$p\equiv 2\pmod{3}$ and $p\equiv \pm 2\pmod{5}$.
Then
$$
F_m\equiv 1\pmod{2}
\quad \text{and} \quad
F_{F_m+2}-1\equiv 0\pmod{F_m}.
$$
\end{theorem}
Since $n\in\cA$ if and only if $F_{n+2}\equiv1\pmod n$ (see \eqref{eq:cong}), and there are infinitely many primes $p$ as in Theorem \ref{thm:FC} by Dirichlet's theorem, there are infinitely many $m$ such that $F_m\in \cA$.
In the following theorem, we obtain necessary conditions for which $F_m\in \cA$.
\begin{theorem}\label{thm:FmA}\noindent
\begin{enumerate}

\item[\upshape(1)] Let $m\geq5$ be odd. If $F_m\in\cA$, then
$$F_m\equiv
\begin{cases}
-3\pmod{4m},&3\nmid m,\\[1mm]
2m-4\pmod{4m},&3\mid m.
\end{cases}
$$
In particular, $m\mid F_m+3$ if $3\nmid m$, and $m\mid F_m+4$ if $3\mid m$.
\item[\upshape(2)] Let $m=p^k\geq5$, where $p$ is an odd prime and $k\geq1$. Then $F_m\notin\cA$.
\item[\upshape(3)] Let $m\geq4$ be even. Then $F_m\in\cA$ if and only if
$$F_m\equiv 0,\ -1,\ \text{or }-3\pmod{2m}.$$
Moreover, $F_m\equiv0\pmod{2m}$ if and only if $12\mid m$ and $m\mid F_m$, and if $3\mid m$, only the residue $0$ is possible.
\item[\upshape(4)] The set $\cA$ contains infinitely many Fibonacci numbers. More precisely:
\begin{enumerate}
\item[\upshape(a)] if $12\mid m$ and $m\mid F_m$, then $F_m\in\cA$, in particular $F_{12\cdot 3^a\cdot 5^b}\in\cA$ for all $a,b\geq0$,
\item[\upshape(b)] if $p$ is a prime with $p\equiv 17$ or $23\pmod{30}$, then $F_{2p}\in\cA$ and $F_{4p}\in\cA$.
\end{enumerate}
\end{enumerate}
\end{theorem}
For example, $F_{12},F_{24},F_{34},F_{36},F_{46},F_{48},F_{60},F_{68},F_{72}\in\cA$, and these are all $F_m\in\cA$ with $3\leq m<90$ apart from $F_3=2$. On the other hand, a computer search shows that there is no odd $m$ with $5\leq m\leq 2\cdot 10^6$ (see \cite{Codes} for Python codes) satisfying the necessary condition of Theorem \ref{thm:FmA}(1). This computer search leads to the following conjecture, to which Theorem \ref{thm:FmA}(2) gives a partial answer.
\begin{conjecture}
If $m\geq 5$ is odd, then $F_m\notin\cA$.
\end{conjecture}

Szemer\'edi \cite{Sz} proved that every set of positive integers with positive upper density contains arbitrarily long arithmetic progressions. The converse is false: by the Green--Tao theorem \cite{GT}, the primes contain arbitrarily long arithmetic progressions, although they have density zero. The set $\cA$ behaves like the primes in this respect, even though $2$ is the only prime in $\cA$ (Theorem \ref{thm:Aeven}(2)): it contains arbitrarily long arithmetic progressions (Theorem \ref{thm:Api}(3)), yet it has density zero. The latter follows from the next result, which concerns all the sums $\sum_{i=1}^nF_{r+di}$ considered in Corollary \ref{Theorem 1} and Theorem \ref{Theorem 2}. 

\begin{theorem}\label{thm:density}
Let $r\geq0$ and $d\geq1$ be integers, and let
$$\cB_{r,d}=\Bigl\{n\in\mathbb{N}:\ n\mid\sum_{i=0}^{n-1}F_{r+id}\Bigr\}.$$
Then $\cB_{r,d}$ has natural density zero. More precisely, for $x\geq3$,
$$\#\{n\leq x:\ n\in\cB_{r,d}\}\ll_d\frac{x}{\sqrt{\log x}},$$
where the implied constant depends only on $d$.
\end{theorem}

\begin{corollary}\label{cor:density}
\noindent
\begin{enumerate}
\item[\upshape(1)] The set $\cA$ has natural density zero. More precisely, $\#\{n\leq x:\ n\in\cA\}\ll x/\sqrt{\log x}$ for $x\geq3$.
\item[\upshape(2)] Let $d\geq1$ and $r\geq-d$ be integers. Then the set of positive integers $n$ with $\displaystyle{n\mid\sum_{i=1}^{n}F_{r+di}}$ has natural density zero. In particular, each of the infinite sets of integers $n$ in Corollary \ref{Theorem 1} and Theorem \ref{Theorem 2} has natural density zero.
\end{enumerate}
\end{corollary}

The paper is organized as follows. In Section \ref{sec:prelim}, we fix notation and collect the known results and auxiliary lemmas needed later. Section \ref{sec:proofs} contains the proofs of the results stated above.

\section{Notations and Preliminaries}\label{sec:prelim}
Throughout, $\mathbb{N}=\{1,2,3,\ldots\}$, and the letter $p$ always denotes a prime number; in particular, sums and products over $p$ run over primes. For a prime $p$ and a nonzero rational number $x$, $\nu_p(x)$ denotes the $p$-adic valuation of $x$, and $\left(\frac{\cdot}{p}\right)$ denotes the Legendre symbol. For a set $\mathcal{M}\subseteq\mathbb{N}$, the \emph{upper density} of $\mathcal{M}$ is
$$\overline{d}(\mathcal{M})=\limsup_{x\to\infty}\frac{\#\{n\leq x:\ n\in\mathcal{M}\}}{x},$$
and $\mathcal{M}$ has \emph{natural density zero} if $\overline{d}(\mathcal{M})=0$, that is, if $\#\{n\leq x:\ n\in\mathcal{M}\}/x\to0$ as $x\to\infty$. For a function $f$ and a positive function $g$, the notations $f\ll g$, $g\gg f$ and $f=O(g)$ all mean that $|f|\leq Cg$ for some constant $C>0$ (in the range of the variables under consideration); a subscript, as in $\ll_d$ or $O_d$, indicates that $C$ may depend on $d$.

Let $A$ and $B$ be two matrices of the same order with integer entries. Then by $A\equiv 0 \pmod n$ we mean that each entry of the matrix is divisible by $n$, and by $A\equiv B \pmod n$, we mean $A-B \equiv 0 \pmod n$. Throughout, $I$ denotes the $2\times2$ identity matrix,
$$
Q=\begin{pmatrix}1&1\\1&0\end{pmatrix},
\qquad
Q^{m}=\begin{pmatrix}F_{m+1}&F_m\\F_m&F_{m-1}\end{pmatrix},
$$
where $m\in\mathbb{Z}$ and $F_{-m}=(-1)^{m+1}F_m$. Note that $L_m=\tr(Q^m)=F_{m+1}+F_{m-1}$, where $\tr$ denotes the trace.

We recall the important facts that for any two positive integers $n,m$, we have $F_n\mid F_m$ if and only if $n\mid m$ (for $n\geq3$), and $\gcd(F_n,F_m)=F_{\gcd(n,m)}$ (see \cite[pp. 197--198]{r13}). We also require the notion of the Pisano period of $n$, which is denoted by $\pi(n)$, it is the smallest positive integer $T$ such that
$$F_{m+T}\equiv F_m \pmod{n}$$
for all integers $m\geq0$, equivalently, $\pi(n)$ is the order of $Q$ modulo $n$. In particular $L_{m+\pi(n)}\equiv L_m\pmod n$ for all $m$. For example, we have $\pi(3)=8$, $\pi(9)=24$, $\pi(5)=20$, $\pi(7)=16$, $\pi(11)=10$, $\pi(2^k)=3\cdot 2^{k-1}$, $\pi(5^k)=4\cdot 5^k$, where $k\geq1$ is an integer, and $\pi(\mathrm{lcm}(a,b))=\mathrm{lcm}(\pi(a),\pi(b))$. We also recall the following properties:
\begin{enumerate}[(1)]
    \item For any prime $p$ and any integers $k\geq1$, $j\geq 0$, we have $\pi(p^{k+j})\mid p^{j}\pi(p^k)$.
    \item For any prime $p\equiv \pm 3 \pmod{10}$, we have $\pi(p)\mid 2(p+1)$.
    \item For any prime $p\equiv \pm 1 \pmod{10}$, we have $\pi(p)\mid p-1$.
\end{enumerate}
Throughout, $\alpha=(1+\sqrt5)/2$ denotes the golden ratio, so that $F_m\leq\alpha^{m-1}$ for all $m\geq1$.

The following results will play a pivotal role in the subsequent sections.
\begin{lemma}[\cite{r4}]\label{lemma1.1}
 Let $m$ be a positive integer greater than $1$. Then $\pi(m)=m$ if and only if $m=24\cdot5^{\lambda-1}$ for some integer $\lambda\geq 1$.
\end{lemma}
\begin{lemma}[{\cite{r6}}]\label{lemma1.7}
For all $n \ge 1$, we have $\nu_5(F_n)=\nu_5(n)$. On the other hand,
$$
\nu_2(F_n)=
\begin{cases}
0, & \text{if } n \equiv 1,2 \pmod{3},\\[4pt]
1, & \text{if } n \equiv 3 \pmod{6},\\[4pt]
3, & \text{if } n \equiv 6 \pmod{12},\\[4pt]
\nu_2(n)+2, & \text{if } n \equiv 0 \pmod{12},
\end{cases}
$$
and
$$
\nu_3(F_n)=
\begin{cases}
0, & \text{if } n \not\equiv 0 \pmod{4},\\[4pt]
\nu_3(n)+1, & \text{if } n \equiv 0 \pmod{4}.
\end{cases}
$$
\end{lemma}
\begin{lemma}[\cite{r7,r8}]\label{lemma1.8}
  If $p\neq 5$ is a prime, then
  $$F_{p-\left(\frac{p}{5}\right)}\equiv0\pmod{p},\quad F_p\equiv \left(\frac{p}{5}\right)\pmod{p}.$$
  Here $\left(\frac{p}{5}\right)=1$ if $p\equiv\pm1\pmod 5$ and $\left(\frac{p}{5}\right)=-1$ if $p\equiv\pm2\pmod 5$.
\end{lemma}

\begin{theorem}[\cite{r5}]\label{thm1.2} We have $m \mid F_n$ and $m \mid F_{n+1}-1$ if and only if $\pi(m)\mid n$.
\end{theorem}
\begin{theorem}[\cite{r14}]\label{thm2.5}
If $F_n=x^2$ for some integer $x$, then $n=0,1,2,12$.
\end{theorem}

We recall the results of Luca and Marques as we shall use them in our subsequent proofs.
A primitive divisor of $F_n$ is a prime factor of $F_n$ that does not divide the product
$$\prod_{j=1}^{n-1} F_j.$$
It is known that such a divisor exists for all $n\neq 1,2,6,12$ (see \cite{r9}). Following Luca and Marques, let us denote by $P(F_n)$ the largest primitive divisor of $F_n$. A positive integer $\lambda$ is called admissible if it satisfies the
following properties:
\begin{enumerate}[(1)]
  \item $\lambda \equiv 0 \pmod{12}$ is a $3$-smooth number (i.e., its largest prime factor is $3$), and both $\lambda-1$ and $\lambda+1$ are prime numbers.
    \item $\lambda(\lambda+1)(\lambda-1)$ divides $F_\lambda$.
    \item $P(F_\lambda)>\lambda+1$.
\end{enumerate}

\begin{lemma}[\cite{r3}]\label{lem2.1}
Let $\lambda$ be admissible. Then there exists a sequence of positive integers
$(k_i)_{i\ge0}$ such that, for each $k=4k_i$, we have
$\gcd\bigl(k,\lambda(\lambda-1)(\lambda+1)\bigr)=4,$
and moreover,
$$\lambda k \mid \frac{F_{\lambda k}}{2F_k}.$$
\end{lemma}

\begin{lemma}[\cite{r3}]\label{lem2.2}
Let $(k_i)_{i\ge0}$ be the sequence provided by Lemma~\ref{lem2.1}.
For each $k=4k_i$, define
$$m:=\frac{F_{\lambda k}}{F_k}.$$
Then both $A_F(m-1)$ and $A_F(m)$ are integers.
\end{lemma}
\begin{remark}\label{remark1}
 Under the notations in Lemmas \ref{lem2.1} and \ref{lem2.2}, during the proof of Conjecture \ref{Conjecture1}, Luca and Marques proved and stated the following:
 \begin{enumerate}
     \item[\upshape(a)] $\dfrac{F_{k(\lambda-1)}F_{k(\lambda+1)}}{F_k}\mid F_{m/2},$
     \item[\upshape(b)] $F_{m+1}\equiv 1 \pmod{m-1},$
     \item[\upshape(c)] $m\mid F_{m/2}$ and $F_{m+2}\equiv 1 \pmod{m}$,
     \item[\upshape(d)] $2^6\mid m$.
 \end{enumerate}
In particular, $m$ is even and $m-1$ is odd, and different values of $k$ give different values of $m$, hence the integers $m$ (respectively $m-1$) form an infinite set of even (respectively odd) positive integers.
\end{remark}

\begin{lemma}\label{Lemma luca}
Under the notations in Lemmas \ref{lem2.1} and \ref{lem2.2}, we have
\begin{enumerate}
    \item[\upshape(1)] $m-1\mid F_{m}, \quad F_{m-1}\equiv 1 \pmod{m-1},$
    \item[\upshape(2)] $m+1\mid F_m$.
\end{enumerate}
\end{lemma}
\begin{proof}
We use the identities $F_a-F_b=F_{(a-b)/2}L_{(a+b)/2}$ and $F_a+F_b=F_{(a+b)/2}L_{(a-b)/2}$, valid whenever $a\equiv b\pmod 2$ and $(a-b)/2$ is even. Since $4\mid k$, the number $k(\lambda-1)/2$ is even.
\begin{enumerate}
    \item[\upshape(1)]
    Set $n=m-1$, that is
$$n=\frac{F_{\lambda k}}{F_k}-1=\frac{F_{\lambda k}-F_k}{F_k}=\frac{F_{k(\lambda-1)/2}L_{k(\lambda+1)/2}}{F_k}.$$
Since $F_{2r}=F_rL_r$, we have
$$\frac{F_{k(\lambda-1)}F_{k(\lambda+1)}}{F_k}\Big/\frac{F_{k(\lambda-1)/2}L_{k(\lambda+1)/2}}{F_k}=L_{k(\lambda-1)/2}F_{k(\lambda+1)/2}\in\mathbb{N},$$
that is, $n\mid \dfrac{F_{k(\lambda-1)}F_{k(\lambda+1)}}{F_k}$.
By Remark \ref{remark1}(a), we have
$$\frac{F_{k(\lambda-1)}F_{k(\lambda+1)}}{F_k}\mid F_{m/2},$$
from which we deduce that $n\mid F_{m}$, that is, $n\mid F_{n+1}$, as $F_{m/2}\mid F_m$. Now by Remark \ref{remark1}(b), we have $F_{n+2}\equiv 1 \pmod{n}$, thus $F_n=F_{n+2}-F_{n+1}\equiv 1 \pmod{n}$.

 \item[\upshape(2)] Set $n=m$. Then
$$n+1=\frac{F_{\lambda k}}{F_k}+1=\frac{F_{\lambda k}+F_k}{F_{k}}=\frac{F_{k(\lambda+1)/2}L_{k(\lambda-1)/2}}{F_k}.$$
Since $F_{2r}=F_rL_r$, we have
$$\dfrac{\displaystyle{\frac{F_{k(\lambda+1)}F_{k(\lambda-1)}}{F_k}}}{n+1}=\frac{F_{k(\lambda +1)}F_{k(\lambda-1)}}{F_{k(\lambda+1)/2}L_{k(\lambda-1)/2}}=L_{k(\lambda+1)/2}F_{k(\lambda-1)/2}\in \mathbb{N},$$
that is
$\displaystyle{n+1}\mid \displaystyle{\frac{F_{k(\lambda+1)}F_{k(\lambda-1)}}{F_k}}.$
By Remark \ref{remark1}(a), we have
$$\frac{F_{k(\lambda+1)}F_{k(\lambda-1)}}{F_k}\mid F_{n/2},$$
and therefore $n+1\mid F_n$, as $F_{n/2}\mid F_n$.
\end{enumerate}
\end{proof}
\begin{proposition}[{\cite[pp. 86--89]{r13}}]\label{Iden}\noindent
    \begin{enumerate}
        \item[\upshape(1)] $$\sum_{i=1}^n F_{3i}=\frac{F_{3n+3}+F_{3n}-2}{4}=\frac{F_{3n+2}-1}{2},$$
        \item[\upshape(2)] $$\sum_{i=1}^nF_{9i}=\frac{F_{9n+9}+F_{9n}-F_9}{L_9}=\frac{F_{9n+9}+F_{9n}-34}{76},$$
        \item[\upshape(3)]  $$\sum_{i=1}^nF_{2i}^2=\frac{3F_{2n+1}^2+2F_{2n+2}^2-6F_{2n}F_{2n+2}-2n-5}{5},$$
        \item[\upshape(4)] $$\sum_{i=1}^nF_{2i-1}^2=\frac{3F_{2n}^2+2F_{2n-1}^2-4F_{2n-2}F_{2n}+2n-2}{5},$$

        \item[\upshape(5)] $$\sum_{i=1}^nF_i^3=\frac{F_{n+2}^3-3F_{n+1}^3+3(-1)^nF_n+2}{4},$$

           \item[\upshape(6)] $$\sum_{i=1}^nF_{2i-1}^3=\frac{F_{2n}^3+3F_{2n}}{4}=\frac{F_{2n}(F_{2n}^2+3)}{4}.$$
    \end{enumerate}
\end{proposition}
\begin{proposition}[{\cite[Ch.~5]{r13}}]\label{prop:basicid}
For all integers $m$ and $n$, we have
\begin{enumerate}
    \item[\upshape(1)] $5F_n^2=L_{2n}-2(-1)^n$,
    \item[\upshape(2)] $5F_n^3=F_{3n}-3(-1)^nF_n$,
    \item[\upshape(3)] $L_{m+n}-L_{m-n}=L_mL_n$ if $n$ is odd,
    \item[\upshape(4)] $L_n^2=5F_n^2+4(-1)^n$,
    \item[\upshape(5)] $2L_{n+1}=L_n+5F_n$,
    \item[\upshape(6)] $F_{2n}=F_nL_n$.
\end{enumerate}
\end{proposition}
\begin{proposition}\label{proposition1}
 For $n\geq1$ we have
  $$\sum_{i=1}^n F_{3i}=\frac{F_{3n+2}-1}{2},\quad \sum_{i=1}^nF_{3i-1}=\frac{F_{3n+1}-1}{2}, \quad \sum_{i=1}^nF_{3i-2}=\frac{F_{3n}}{2}.$$
\end{proposition}
\begin{proof}
 From Proposition \ref{Iden}, we find that
   $$\sum_{i=1}^n F_{3i}=\frac{F_{3n+3}+F_{3n}-2}{4}=\frac{F_{3n+2}-1}{2}.$$
The other two identities follow by induction on $n$.
\end{proof}
\begin{proposition}\label{proposition2}
For $n\geq1$, we have
 $$\sum_{i=1}^{n} F_{3i}^{2}=\frac{L_{6n+3}-4(-1)^n}{20},$$
 $$\sum_{i=1}^{n} F_{3i-1}^{2}=\frac{L_{6n+1}+4(-1)^n-5}{20}, \quad \sum_{i=1}^{n} F_{3i-2}^{2}=\frac{L_{6n-1}+5-4(-1)^n}{20}.$$
\end{proposition}
\begin{proof}
     From Proposition \ref{prop:basicid}(1), we find that
     $5F_n^2=L_{2n}-2(-1)^n.$
     Therefore,
     $$\sum_{i=1}^nF_{3i}^2=\frac{1}{5}\sum_{i=1}^n\bigl(L_{6i}-2(-1)^{3i}\bigr)=\frac{1}{5}\sum_{i=1}^nL_{6i}-\frac{(-1)^n-1}{5}.$$
  From Proposition \ref{prop:basicid}(3), we find that for odd $n$ the identity $L_{m+n}-L_{m-n}=L_mL_n$ holds, setting $m=6i$ and $n=3$, we get
  $$4L_{6i}=L_{6i+3}-L_{6i-3},$$
  which gives
  $$\sum_{i=1}^nL_{6i}=\frac{1}{4}\sum_{i=1}^n(L_{6i+3}-L_{6i-3})=\frac{L_{6n+3}-4}{4},$$
  consequently, we have
  $$\sum_{i=1}^nF_{3i}^2=\frac{L_{6n+3}-4}{20}-\frac{(-1)^n-1}{5}=\frac{L_{6n+3}-4(-1)^n}{20}.$$
  Similarly, we get
  $$\sum_{i=1}^nF_{3i-1}^2=\frac{1}{5}\sum_{i=1}^n\bigl(L_{6i-2}-2(-1)^{3i-1}\bigr)=\frac{1}{5}\sum_{i=1}^nL_{6i-2}-\frac{1-(-1)^{n}}{5}.$$
Now, setting  $m=6i-2$ and $n=3$, we get
$4L_{6i-2}=L_{6i+1}-L_{6i-5},$
which gives
$$\sum_{i=1}^nL_{6i-2}=\frac{1}{4}\sum_{i=1}^n(L_{6i+1}-L_{6i-5})=\frac{L_{6n+1}-L_1}{4}=\frac{L_{6n+1}-1}{4},$$
consequently, we have
  $$\sum_{i=1}^nF_{3i-1}^2=\frac{L_{6n+1}-1}{20}-\frac{1-(-1)^n}{5}=\frac{L_{6n+1}+4(-1)^n-5}{20}.$$
  The third identity is proved in the same way, using $4L_{6i-4}=L_{6i-1}-L_{6i-7}$ and $L_{-1}=-1$.
  \end{proof}

\begin{proposition}\label{proposition3}
For $n\geq1$, we have
  $$\sum_{i=1}^{n} F_{3i}^{3}=\frac{1}{5}\left(\frac{F_{9n+9}+F_{9n}-34}{76}\right)-\frac{3}{5}\left(\frac{(-1)^nF_{3n+1}-1}{2}\right),$$
$$\sum_{i=1}^{n} F_{3i-1}^{3}=\frac{1}{5}\left(\frac{10F_{9n+1}+2F_{9n}-10}{76}\right)+\frac{3}{5}\cdot\frac{(-1)^nF_{3n}}{2}, $$
$$\sum_{i=1}^{n} F_{3i-2}^{3}=\frac{1}{5}\left(\frac{14F_{9n}-6F_{9n+1}+6}{76}\right)-\frac{3}{5}\left(
\frac{(-1)^nF_{3n-1}-1}{2}\right).$$
\end{proposition}
\begin{proof}
By Proposition \ref{prop:basicid}(2), we have
$F_{3n}=5F_n^3+3(-1)^nF_n,$
setting $n=3i$, we get
$$F_{3i}^3=\frac{F_{9i}-3(-1)^iF_{3i}}{5}.$$
Therefore,
$$\sum_{i=1}^nF_{3i}^3=\frac{1}{5}\sum_{i=1}^nF_{9i}
-\frac{3}{5}\sum_{i=1}^n(-1)^iF_{3i}.$$
 From Proposition \ref{Iden}, we have
$$\sum_{i=1}^nF_{9i}=\frac{F_{9n+9}+F_{9n}-F_9}{L_9}=\frac{F_{9n+9}+F_{9n}-34}{76}.$$
Since
$F_{3i+1}+F_{3i-2}=2F_{3i},$
we obtain
$$(-1)^iF_{3i}=\frac{1}{2}\left((-1)^iF_{3i+1}-(-1)^{i-1}F_{3i-2}
\right).$$
Therefore, the sum telescopes and
\begin{align*}
\sum_{i=1}^n(-1)^iF_{3i}=\frac{1}{2}\sum_{i=1}^n\left((-1)^iF_{3i+1}
-(-1)^{i-1}F_{3i-2}\right)=\frac{(-1)^nF_{3n+1}-F_1}{2}=
\frac{(-1)^nF_{3n+1}-1}{2}.
\end{align*}
This proves the first identity. The other two identities are obtained in the same way from $F_{3i-1}^3=\bigl(F_{9i-3}+3(-1)^{i}F_{3i-1}\bigr)/5$ and $F_{3i-2}^3=\bigl(F_{9i-6}-3(-1)^{i}F_{3i-2}\bigr)/5$.
\end{proof}

The next lemma, together with its refinement Lemma \ref{lem:lift5} for $p=5$, is the main tool in the proofs of Theorems \ref{thm:lambda}, \ref{Theorem 4} and \ref{Theorem 7}.
\begin{lemma}\label{Lemma 3.10}
Let $s$ and $t$ be positive integers, and $p$ be a prime. Let $A$ be an integer matrix such that
$$A\equiv I\pmod p.$$
If $p^s\mid t$, then
$$\sum_{j=0}^{t-1}A^j\equiv0\pmod{p^s}.$$
\end{lemma}

\begin{proof}
Since $A\equiv I\pmod p$, we have
$A=I+pC$ for some integer matrix $C$.
Hence, by the binomial theorem,
$$A^j=(I+pC)^j=\sum_{\ell=0}^{j}\binom{j}{\ell}p^\ell C^\ell.$$
Consequently, using $\displaystyle{\sum_{j=\ell}^{t-1}\binom{j}{\ell}=\binom{t}{\ell+1}}$,
\begin{align}\label{eq:binsum}
\sum_{j=0}^{t-1}A^j = \sum_{j=0}^{t-1}\sum_{\ell=0}^{j}\binom{j}{\ell}p^\ell C^\ell=
\sum_{\ell=0}^{t-1}p^\ell C^\ell\sum_{j=\ell}^{t-1}\binom{j}{\ell}=
\sum_{\ell=0}^{t-1}
p^\ell
\binom{t}{\ell+1}
C^\ell.
\end{align}
To complete the proof, it is sufficient to prove that
$$p^s\mid p^\ell\binom{t}{\ell+1}
\qquad \ \text{for} \  \ell = 0, 1, \ldots, t-1.$$
Clearly, this holds for  $\ell=0$. Now, let $\ell\geq1$. Since
$$\binom{t}{\ell+1}=\frac{t}{\ell+1}\binom{t-1}{\ell},$$
we have
$$\nu_p\left(p^\ell\binom{t}{\ell+1}
\right)=\ell+\nu_p(t)-\nu_p(\ell+1)+
\nu_p\left(\binom{t-1}{\ell}\right)\geq \nu_p(t)+\bigl(\ell-\nu_p(\ell+1)\bigr).$$
Since $p^{\ell}\geq 2^\ell\geq \ell+1$, we have $\nu_p(\ell+1)\leq \ell$, and $\nu_p(t)\geq s$. Therefore
$\nu_p\left(p^\ell\binom{t}{\ell+1}\right)\geq s,$ which proves the lemma.
\end{proof}

We next prove a refinement of Lemma \ref{Lemma 3.10} for $p=5$.

\begin{lemma}\label{lem:lift5}
Let $A$ be a $2\times2$ integer matrix such that $A^5\equiv I\pmod 5$ and
$$S:=A+A^2+A^3+A^4+A^5\equiv 0\pmod 5.$$
Then for every $k\geq1$,
$$\sum_{i=1}^{5^k}A^i\equiv 5^{k-1}S\pmod{5^{k+1}}.$$
In particular, $\displaystyle{\sum_{i=1}^{5^k}A^i\equiv 0\pmod{5^{k}}}$.
Moreover, if $4\mid e$, then the matrix $A=Q^e$ satisfies both hypotheses.
\end{lemma}
\begin{proof}
Grouping the indices according to their residue modulo $5$, we get
$$\sum_{i=1}^{5^k}A^i=S\sum_{j=0}^{5^{k-1}-1}(A^5)^j.$$
Write $A^5=I+5C$ and $t=5^{k-1}$. By \eqref{eq:binsum}, $\displaystyle{\sum_{j=0}^{t-1}(A^5)^j=\sum_{\ell\geq0}5^\ell\binom{t}{\ell+1}C^\ell}$. The term with $\ell=0$ is $tI=5^{k-1}I$, and, as in the proof of Lemma \ref{Lemma 3.10}, for $\ell\geq1$ we have $\nu_5\bigl(5^\ell\binom{t}{\ell+1}\bigr)\geq (k-1)+\ell-\nu_5(\ell+1)\geq k$, because $5^{\ell-1}\geq \ell+1$ for $\ell\geq 2$ and $\nu_5(2)=0$. Hence $\displaystyle{\sum_{j=0}^{t-1}(A^5)^j=5^{k-1}I+5^kE}$ for some integer matrix $E$, and
$$\sum_{i=1}^{5^k}A^i=5^{k-1}S+5^kSE\equiv 5^{k-1}S\pmod{5^{k+1}},$$
since $S\equiv0\pmod 5$.

Now let $e=4v$. A direct computation gives $Q^4=\begin{pmatrix}5&3\\3&2\end{pmatrix}\equiv I+N\pmod 5$ with $N=\begin{pmatrix}4&3\\3&1\end{pmatrix}$ and $N^2\equiv0\pmod 5$. Hence $A=Q^{4v}\equiv (I+N)^v\equiv I+vN\pmod 5$ and, more generally, $A^j\equiv I+jvN\pmod 5$ for all $j\geq0$. Therefore $A^5\equiv I+5vN\equiv I\pmod 5$ and
$$S\equiv 5I+(1+2+3+4+5)vN=5I+15vN\equiv 0\pmod 5.$$
\end{proof}

Combining Lemma \ref{Lemma 3.10} with $\pi(2)=3$ and $\pi(3)=8$, we obtain the following congruence, which produces the even integers $n$ in Theorems \ref{thm:lambda}, \ref{Theorem 4} and \ref{Theorem 7}.
\begin{lemma}\label{lem:key}
Let $r\geq0$ and $d\geq 1$ be integers, let $\lambda\in\{-1,1\}$, and let $k\geq1$. Then
$$\sum_{i=1}^{8\cdot 3^{k}}\lambda^{i}Q^{r+di}\equiv 0\pmod{8\cdot 3^k}.$$
\end{lemma}
\begin{proof}
Set $A:=\lambda Q^{d}$, so that $\lambda^{i}Q^{r+di}=Q^{r}A^{i}$ for every $i\ge 1$.
Since $\lambda^{8}=1$ and $Q^{8}\equiv I\pmod 3$ (as $\pi(3)=8$), we have
\[
A^{8}=(\lambda Q^{d})^{8}=\lambda^{8}Q^{8d}=(Q^{8})^{d}\equiv I\pmod 3.
\]
Grouping the indices according to the residue of the exponent modulo $8$, we obtain
\begin{align*}
\sum_{i=1}^{8\cdot 3^{k}}\lambda^{i}Q^{r+di}
=Q^{r}\sum_{i=1}^{8\cdot 3^{k}}A^{i}
=Q^{r}\sum_{j=0}^{3^{k}-1}\sum_{t=1}^{8}A^{8j+t}
=Q^{r}\left(\sum_{t=1}^{8}A^{t}\right)\left(\sum_{j=0}^{3^{k}-1}(A^{8})^{j}\right).
\end{align*}
By Lemma~\ref{Lemma 3.10}, $\displaystyle{\sum_{j=0}^{3^{k}-1}(A^{8})^{j}\equiv 0\pmod{3^{k}}}$, hence
 $\sum_{i=1}^{8\cdot 3^{k}}\lambda^{i}Q^{r+di}\equiv 0\pmod{3^{k}}$.

Since $-1\equiv 1\pmod 2$, we have $A=\lambda Q^{d}\equiv Q^{d}\pmod 2$. As $Q^{3}\equiv I\pmod 2$ (since $\pi(2)=3$), it follows that
$A^{3}\equiv (Q^{3})^{d}\equiv I\pmod 2.$
Grouping the indices according to the residue of the exponent modulo $3$, we obtain
\begin{align*}
\sum_{i=1}^{8\cdot 3^{k}}\lambda^{i}Q^{r+di}
=Q^{r}\sum_{j=0}^{8\cdot 3^{k-1}-1}\left(A^{3j+1}+A^{3j+2}+A^{3j+3}\right)
=Q^{r}\left(A+A^{2}+A^{3}\right)\sum_{j=0}^{8\cdot 3^{k-1}-1}(A^{3})^{j}.
\end{align*}
By Lemma~\ref{Lemma 3.10} (with $p=2$, $s=3$), $\displaystyle{\sum_{j=0}^{8\cdot 3^{k-1}-1}(A^{3})^{j}\equiv 0\pmod 8}$, hence
 $\displaystyle{\sum_{i=1}^{8\cdot 3^{k}}\lambda^{i}Q^{r+di}\equiv 0\pmod 8}$. Since $\gcd(8,3^{k})=1$, the lemma follows.
\end{proof}

To prove Theorems \ref{thm:FC} and \ref{thm:FmA}, we need the following auxiliary lemma.

\begin{lemma}\label{lem:Fk1}
\noindent
\begin{enumerate}
\item[\upshape(1)] Let $m\geq4$ be even. Then $F_k\equiv1\pmod{F_m}$ if and only if $k\equiv1,\ 2,\ -1\pmod{2m}$.
\item[\upshape(2)] Let $m\geq5$ be odd. Then $F_k\equiv1\pmod{F_m}$ if and only if $k\equiv1,\ 2,\ 2m-2,\ -1\pmod{4m}$.
\end{enumerate}
\end{lemma}

\begin{proof}
Put $M=F_m$ and $c=F_{m-1}$. Since $F_{m+1}=F_m+F_{m-1}\equiv c\pmod M$, we have
$$Q^m=
\begin{pmatrix}
F_{m+1}&F_m\\
F_m&F_{m-1}
\end{pmatrix}
\equiv
cI
\pmod M.$$
Cassini's identity $F_{m+1}F_{m-1}-F_m^2=(-1)^m$ gives $c^2\equiv(-1)^m\pmod M$. In particular $c$ is invertible modulo $M$, $Q^{2m}\equiv I\pmod M$ if $m$ is even, and $Q^{2m}\equiv -I$, $Q^{4m}\equiv I\pmod M$ if $m$ is odd. Hence $F_k\bmod M$ depends only on $k$ modulo $2m$ (respectively $4m$). Every integer $k$ can be written as $k=qm+s$ with $0\leq s<m$, and then $Q^k\equiv c^qQ^s$, so that $F_k\equiv c^qF_s\pmod M$. Note that $0\leq F_s\leq F_{m-1}<M$ and that $F_s$ determines $s$ when $F_s\geq2$. Let $m$ be even, so $c^2\equiv1$, and take $q\in\{0,1\}$. If $q=0$, then $F_k\equiv1$ means $F_s=1$, i.e., $s\in\{1,2\}$. If $q=1$, then $cF_s\equiv1$ means $F_s\equiv c=F_{m-1}\pmod M$, i.e., $F_s=F_{m-1}$, so $s=m-1$ (note $F_{m-1}\geq 2$) and $k=2m-1$. Hence $F_k\equiv1\pmod M$ if and only if $k\equiv1,2,-1\pmod{2m}$. Let $m$ be odd, so $c^2\equiv-1$, and take $q\in\{0,1,2,3\}$. If $q=0$, then $s\in\{1,2\}$ as before. If $q=2$, then $-F_s\equiv1$, i.e., $F_s=M-1$, which is impossible because $F_s\leq F_{m-1}=M-F_{m-2}<M-1$ (as $F_{m-2}\geq2$). If $q=1$, then $cF_s\equiv1$, multiplying by $-c$ we get $F_s\equiv-c\equiv M-F_{m-1}=F_{m-2}\pmod M$, hence $s=m-2$ and $k=2m-2$. If $q=3$, then $c^3=c\cdot c^2\equiv -c$, so $-cF_s\equiv 1$, multiplying by $c$ we get $F_s\equiv c=F_{m-1}$, hence $s=m-1$ and $k=4m-1$. This proves (2).
\end{proof}

The following lemma will be used in the proofs of Theorems \ref{thm:Aeven}(2) and \ref{thm:FmA}(2).
\begin{lemma}\label{lem:Fpk}
Let $p\neq 2,5$ be a prime, let $\varepsilon=\left(\frac p5\right)$ and $k\geq1$. Then
$$Q^{p^k}\equiv\begin{cases}Q, & \text{if } \varepsilon^k=1,\\ -Q^{-1}, & \text{if } \varepsilon^k=-1,\end{cases}\pmod p.$$ Moreover,  $F_{p^k}\equiv\varepsilon^k\pmod p$.
\end{lemma}

\begin{proof}
If $\varepsilon=1$, Lemma \ref{lemma1.8} gives $F_{p-1}\equiv0$ and $F_p\equiv1\pmod p$, hence $F_{p-2}=F_p-F_{p-1}\equiv1$ and
$$Q^{p-1}=\left(\begin{matrix}F_p&F_{p-1}\\F_{p-1}&F_{p-2}\end{matrix}\right)\equiv I\pmod p.$$
Thus $Q^p\equiv Q$, and by induction $Q^{p^k}=(Q^{p^{k-1}})^p\equiv Q^p\equiv Q\pmod p$.

If $\varepsilon=-1$, Lemma \ref{lemma1.8} gives $F_{p+1}\equiv0$ and $F_p\equiv-1\pmod p$, hence $F_{p+2}=F_{p+1}+F_p\equiv-1$ and
$$Q^{p+1}=\left(\begin{matrix}F_{p+2}&F_{p+1}\\F_{p+1}&F_{p}\end{matrix}\right)\equiv -I\pmod p,$$ i.e., $Q^p\equiv-Q^{-1}$. Since $p$ is odd,
$$\bigl(-Q^{-1}\bigr)^p=-\bigl(Q^{p}\bigr)^{-1}\equiv-\bigl(-Q^{-1}\bigr)^{-1}=Q\pmod p .$$
Hence $Q^{p^k}\equiv -Q^{-1}$ for odd $k$ and $Q^{p^k}\equiv Q$ for even $k$, by induction on $k$. Finally, the $(1,2)$-entries of $Q$ and $-Q^{-1}=\left(\begin{smallmatrix}0&-1\\-1&1\end{smallmatrix}\right)$ are $1$ and $-1$, which gives $F_{p^k}\equiv\varepsilon^k\pmod p$.
\end{proof}

For the proof of Theorem \ref{thm:density} we work in a finite field. Let $p\neq5$ be a prime. The polynomial $X^2-X-1$ has discriminant $5\not\equiv0\pmod p$, so it has two distinct roots $\theta,\theta'$ in the field $\mathbb{F}_{p^2}$, with $\theta+\theta'=1$, $\theta\theta'=-1$ and $\theta-\theta'\neq0$. Both sides of
\begin{equation}\label{eq:binetp}
F_m\equiv\frac{\theta^m-\theta'^m}{\theta-\theta'}\quad \text{where }m\geq0
\end{equation}
satisfy the Fibonacci recurrence and agree for $m=0,1$, so \eqref{eq:binetp} holds in $\mathbb{F}_{p^2}$.

\begin{lemma}\label{lem:fieldp}
Let $p\neq5$ be a prime. Then
\begin{enumerate}
\item[\upshape(1)] for every integer $T\geq1$, we have $\theta^T=\theta'^T=1$ if and only if $\pi(p)\mid T$,
\item[\upshape(2)] $\gcd(p,\pi(p))=1$,
\item[\upshape(3)] $\pi(p)\geq\log p/\log\alpha$.
\end{enumerate}
\end{lemma}

\begin{proof}
(1) By Theorem \ref{thm1.2}, $\pi(p)\mid T$ if and only if $F_T\equiv0\pmod p$ and $F_{T+1}\equiv1\pmod p$. If $\theta^T=\theta'^T=1$, then \eqref{eq:binetp} gives $F_T\equiv0\pmod p$ and $F_{T+1}\equiv1\pmod p$. Conversely, $F_T\equiv0\pmod p$ gives $\theta^T=\theta'^T$, and then $F_{T+1}\equiv1\pmod p$ gives $\theta^T(\theta-\theta')=\theta-\theta'$, so $\theta^T=\theta'^T=1$.

(2) We have $\pi(2)=3$. For $p\neq2,5$, the Pisano period $\pi(p)$ divides $p-1$ or $2(p+1)$ by properties (2) and (3) above. This proves $\gcd(p,\pi(p))=1$.

(3) By Theorem \ref{thm1.2}, $p\mid F_{\pi(p)}$, and $F_{\pi(p)}\geq1$. Hence $p\leq F_{\pi(p)}\leq\alpha^{\pi(p)-1}$, which gives $\pi(p)\geq1+\log p/\log\alpha$.
\end{proof}

\begin{lemma}\label{lem:classes}
Let $r\geq0$ and $d\geq1$ be integers, and let $p\neq5$ be a prime with $\pi(p)\nmid d$. Let
$$S(n)=\sum_{i=0}^{n-1}F_{r+id}.$$
Then the set
$$E_p=\{n\geq1:\ p\mid n\ \text{and}\ p\mid S(n)\}$$
is contained in a union of at most $4d$ residue classes modulo $p\,\pi(p)$.
\end{lemma}

\begin{proof}
Set $X_n=\theta^{dn}$ and $Y_n=\theta'^{dn}$. Since $\theta'=-\theta^{-1}$, we have $Y_n=\varepsilon_nX_n^{-1}$ with $\varepsilon_n=(-1)^{dn}\in\{1,-1\}$. By Lemma \ref{lem:fieldp}(1) and the hypothesis $\pi(p)\nmid d$, at least one of $\theta^d,\theta'^d$ is different from $1$. By \eqref{eq:binetp}, we can write
$$(\theta-\theta')\,S(n)=\sum_{i=0}^{n-1}\theta^{r+id}-\sum_{i=0}^{n-1}\theta'^{r+id}.$$
We claim that for $n\in E_p$ the pair $(X_n,Y_n)$ takes at most $4$ values.

\emph{Case 1: $\theta^d\neq1$ and $\theta'^d\neq1$.} Then
$$(\theta-\theta')\,S(n)=a(X_n-1)-b(Y_n-1),\qquad a=\frac{\theta^r}{\theta^d-1}\neq0,\quad b=\frac{\theta'^r}{\theta'^d-1}\neq0.$$
If $p\mid S(n)$, then $aX_n-b\varepsilon_nX_n^{-1}=a-b$, i.e.,
$$aX_n^2-(a-b)X_n-b\varepsilon_n=0.$$
For each of the two values of $\varepsilon_n$, this quadratic equation has a nonzero leading coefficient, hence at most $2$ roots. So $X_n$ takes at most $4$ values, and $Y_n=\varepsilon_nX_n^{-1}$ is determined by $X_n$ and $\varepsilon_n$. Hence the pair $(X_n,Y_n)$ takes at most $4$ values.

\emph{Case 2: exactly one of $\theta^d,\theta'^d$ equals $1$}, say $\theta^d=1\neq\theta'^d$ (the other case is symmetric). Then $\sum_{i=0}^{n-1}\theta^{r+id}=n\theta^r$, which vanishes in $\mathbb{F}_{p^2}$ because $p\mid n$. Hence $(\theta-\theta')S(n)=-b(Y_n-1)$ with $b\neq0$ as above, and $p\mid S(n)$ forces $Y_n=1$. Since also $X_n=(\theta^d)^n=1$, the pair $(X_n,Y_n)$ takes only the value $(1,1)$.

This proves the claim. Now, by Lemma \ref{lem:fieldp}(1), $(X_n,Y_n)=(X_{n'},Y_{n'})$ if and only if $\pi(p)\mid d(n-n')$, i.e., if and only if $n\equiv n'\pmod{\pi(p)/g}$ with $g=\gcd(d,\pi(p))$. Hence each value of the pair corresponds to exactly $g\leq d$ residue classes of $n$ modulo $\pi(p)$, and $E_p$ is contained in at most $4d$ residue classes of $n$ modulo $\pi(p)$, intersected with the condition $n\equiv0\pmod p$. Since $\gcd(p,\pi(p))=1$ by Lemma \ref{lem:fieldp}(2), the Chinese Remainder Theorem shows that $E_p$ is contained in at most $4d$ residue classes modulo $p\,\pi(p)$.
\end{proof}

\begin{lemma}\label{lem:sumpi}
For $y\geq3$, we have
$$\sum_{p>y}\frac{1}{p\,\pi(p)}=O\left(\frac{1}{\log y}\right).$$
\end{lemma}

\begin{proof}
By Lemma \ref{lem:fieldp}(3) (and $\pi(5)=20$), $\displaystyle{\sum_{p>y}\frac{1}{p\,\pi(p)}}\leq \log\alpha\displaystyle{\sum_{p>y}\frac{1}{p\log p}}$.
Let $\Pi(t)$ denote the number of primes not exceeding $t$ (we avoid $\pi(t)$, which denotes the Pisano period). By Abel's summation formula, applied to $a_n=1$ if $n$ is prime and $a_n=0$ otherwise (so that $\sum_{n\leq t}a_n=\Pi(t)$) and to $f(t)=\dfrac{1}{t\log t}$, we have
$$
\sum_{y<p\le x}\frac{1}{p\log p}
=\frac{\Pi(x)}{x\log x}-\frac{\Pi(y)}{y\log y}
+\int_y^x \Pi(t)\,\frac{\log t+1}{t^{2}\log^{2}t}\,dt .
$$
Since $\Pi(t)=O( t/\log t)$, the first term tends to $0$ as $x\to\infty$ and the second is $O(1/\log^2 y)$ and since $\log t+1\le 2\log t$
for $t\ge e$, letting $x\to\infty$ gives
$$
\sum_{p>y}\frac{1}{p\log p}
=O\left( \int_y^\infty \frac{\Pi(t)}{t^{2}\log t}\,dt
\right)=O\left(\int_y^\infty \frac{t/\log t}{t^{2}\log t}\,dt\right)
=O\left(\int_y^\infty \frac{dt}{t\log^{2}t}\right)
=O\left(\frac{1}{\log y}\right).
$$
\end{proof}

Now we establish some auxiliary results on the set
\[
\cA=\left\{n\in\mathbb{N}\setminus\{1\}:
n\mid\sum_{i=1}^{n}F_i\right\}=\{2,24,48,72,77,96,120,144,192,216,240,\ldots\}.
\]
Since
\[
\sum_{i=1}^nF_i=F_{n+2}-1,
\]
we have
\begin{equation}\label{eq:cong}
n\in\cA \quad \text{if and only if}\quad  F_{n+2}\equiv1\pmod n.
\end{equation}
We shall also use that the Fibonacci sequence modulo $3$, $4$ and $5$ is periodic with period $8$, $6$ and $20$, respectively, with
$$(F_k\bmod 3)_{k=0}^{7}=(0,1,1,2,0,2,2,1),\qquad (F_k\bmod 4)_{k=0}^{5}=(0,1,1,2,3,1),$$
so that
\begin{align}\label{eq:F1mod345}
\begin{split}
F_k\equiv1\pmod3 \quad \text{if and only if}\quad k\equiv1,2,7\pmod8,\\
F_k\equiv1\pmod4\quad \text{if and only if}\quad k\equiv1,2,5\pmod6,\\
F_k\equiv1\pmod5\quad \text{if and only if}\quad k\equiv1,2,8,19\pmod{20}.
\end{split}
\end{align}
Finally, by Lemma \ref{lemma1.7}, $F_k$ is even if and only if $3\mid k$.

\section{Proofs}\label{sec:proofs}


\begin{proof}[Proof of Theorem \ref{Theorem 2}]Note that  \cite[p. 86]{r13}
$$\sum_{i=1}^{n}F_{2i-1}=F_{2n}=F_{n}L_{n},\quad  \sum_{i=1}^{n}F_{2i}=F_{2n+1}-1.$$
For $n=24\cdot 5^{k-1}$, where $k\geq1$, we have $\pi(n)=n$ by Lemma \ref{lemma1.1}. Therefore, by Theorem \ref{thm1.2} we have $n\mid F_{2n}$ and $n\mid F_{2n+1}-1$. This proves that there exist infinitely many even positive integers $n$ such that
$$\frac{1}{n}\sum_{i=1}^{n}F_{2i-1},\quad  \frac{1}{n}\sum_{i=1}^{n}F_{2i}\in \mathbb{N}.$$
Let $p$ be a prime such that $p\equiv\pm2\pmod{5}$. Then $\left(\frac{p}{5}\right)=-1$, and from Lemma \ref{lemma1.8} we get $p\mid F_{p+1}$ and $F_p\equiv -1\pmod p$. Since
$F_{2p+1}=F_{p}^2+F_{p+1}^2$, we obtain
$$F_{2p+1}\equiv F_p^2\equiv 1\pmod{p},$$
that is, $\displaystyle{\frac{1}{p} \sum_{i=1}^pF_{2i}}\in \mathbb{N}$, and there are infinitely many such primes $p$ by Dirichlet's theorem. Set $n=5^{k}$, where $k\geq1$ is an integer. Then from Lemma \ref{lemma1.7} we get $\nu_5(F_n)=\nu_5(n)=k$, i.e., $5^{k}\mid F_n\mid F_{2n}$ for all $k\geq1$. Hence $\displaystyle{\frac{1}{n} \sum_{i=1}^nF_{2i-1}}\in \mathbb{N}$. This proves (1) and (2).\\

  From Proposition \ref{proposition1}, we have
  \begin{align}\label{eq23i}
   \sum_{i=1}^n F_{3i}=\frac{F_{3n+2}-1}{2},\quad \sum_{i=1}^nF_{3i-1}=\frac{F_{3n+1}-1}{2}, \quad \sum_{i=1}^nF_{3i-2}=\frac{F_{3n}}{2}.
  \end{align}

   Set $n=2^{k}$, where $k\geq1$.
   Since $\pi(2^{k+1})=3\cdot 2^{k}=3n$, it follows from Theorem \ref{thm1.2} that $2^{k+1}\mid F_{3n}$ and $F_{3n+1}\equiv 1\pmod{2^{k+1}}$, and thus $F_{3n+2}=F_{3n+1}+F_{3n}\equiv 1 \pmod{2^{k+1}}$. This together with (\ref{eq23i}), we have

$$\frac{1}{n}\sum_{i=1}^n F_{3i},\quad \frac{1}{n}\sum_{i=1}^nF_{3i-1}, \quad \frac{1}{n}\sum_{i=1}^nF_{3i-2}\in \mathbb{N}.$$
From $Q^{3n+1}=(Q^n)^3Q$ (or by Binet's formula), we have
$$F_{3n+1}
=F_{n+1}^3+3F_{n+1}F_n^2-F_n^3,\quad F_{3n+2}=F_{n+1}^3+3F_{n+1}^2F_n+F_n^3.$$
Set $n=m-1$ in Lemma \ref{lem2.2}, then $n$ is odd, and these integers form an infinite set by Remark \ref{remark1}. Since by Lemma \ref{Lemma luca}, we have $n\mid F_{n+1}$ and $F_n\equiv 1 \pmod{n}$, we get
$$ F_{3n+2}-1=F_{n+1}^3+3F_{n+1}^2F_n+F_n^3-1\equiv 0\pmod{n},$$
that is, $n\mid F_{3n+2}-1$, and since $n$ is odd,
$$ \frac{1}{n}\sum_{i=1}^nF_{3i}\in \mathbb{N}.$$
Now take a prime $p$ such that $p\equiv \pm 3 \pmod{10}$ (there are infinitely many such primes by Dirichlet's theorem). By Lemma \ref{lemma1.8}, we have $p\mid F_{p+1}$ and $F_p\equiv -1 \pmod{p}$. Therefore
$$F_{3p+1}-1=F_{p+1}^3+3F_{p+1}F_p^2-F_p^3-1\equiv1-1\equiv0\pmod{p},$$
that is,
$$ \frac{1}{p}\sum_{i=1}^pF_{3i-1}\in \mathbb{N}.$$
Set $n=5^k$, where $k\geq1$. Then by Lemma \ref{lemma1.7}, we have $\nu_5(F_{3n})=\nu_5(3n)=k$, that is, $n\mid F_{3n}$, consequently
$$ \frac{1}{n}\sum_{i=1}^nF_{3i-2}\in \mathbb{N}.$$
This proves (3), (4), and (5).
\end{proof}

\begin{proof}[Proof of Theorem \ref{thm:lambda}]
By Lemma \ref{lem:key}, $\displaystyle{\sum_{i=1}^{8\cdot 3^{k}}\lambda^{i}Q^{r+di}\equiv 0\pmod{8\cdot 3^k}}$. Since the $(1,2)$-entry of $\lambda^iQ^{r+di}$ is $\lambda^{i}F_{r+di}$, extracting the $(1,2)$-entry, we get
\[
8\cdot 3^{k}\mid\sum_{i=1}^{8\cdot 3^{k}}\lambda^{i}F_{r+di}.
\]
This holds for every $k\ge 1$ and each $\lambda\in\{-1,1\}$, hence there exist infinitely many
even positive integers $n$ (namely $n=8\cdot 3^{k}$) such that
$\displaystyle{n\mid\sum_{i=1}^{n}\lambda^{i}F_{r+di}}.$

Assume now that $d=4v$, let $n=5^{k}$ and $A=Q^{d}$. By Lemma \ref{lem:lift5}, $A^5\equiv I\pmod 5$ and $A+A^{2}+A^{3}+A^{4}+A^{5}\equiv 0\pmod 5$, and
\[
\sum_{i=1}^{5^{k}}Q^{r+di}=Q^{r}\sum_{i=1}^{5^{k}}A^{i}\equiv 0\pmod{5^{k}}.
\]
Extracting the $(1,2)$-entry, we get $\displaystyle{5^{k}\mid\sum_{i=1}^{5^{k}}F_{r+di}}$. This completes the proof.
\end{proof}
\begin{proof}[Proof of Corollary \ref{Theorem 1}]
This is the case $\lambda=1$ of Theorem \ref{thm:lambda}.
\end{proof}

\begin{proof}[Proof of Theorem \ref{Theorem 5}] Note that
\begin{equation}\label{Equation-Theorem-1.8}
\sum_{i=1}^{n}F_{i}^2=F_nF_{n+1}.
\end{equation}
Note that $n\mid F_n$ when $n=24\cdot 5^{k-1}$ (by Lemma \ref{lemma1.1} and Theorem \ref{thm1.2}) and when $n=5^k$ (by Lemma \ref{lemma1.7}). Therefore from \eqref{Equation-Theorem-1.8}, there exist  infinitely many even and odd positive integers $n$ such that
$$n\mid \sum_{i=1}^{n}F_{i}^2.$$
This proves (1).

From Proposition \ref{Iden}, we have
    $$\sum_{i=1}^nF_{2i}^2=\frac{3F_{2n+1}^2+2F_{2n+2}^2-6F_{2n}F_{2n+2}-2n-5}{5}.$$
 Set $n=3\cdot 2^{k+1}$ for $k\geq1$. Since $\pi(3)=8\mid 2n$, by Theorem \ref{thm1.2}, we have $3\mid F_{2n}$ and $F_{2n+1}\equiv 1 \pmod{3}$. Therefore,
$F_{2n+2}=F_{2n+1}+F_{2n}\equiv 1 \pmod{3}$. Since $3\mid 2n$,
$$3F_{2n+1}^2+2F_{2n+2}^2-6F_{2n}F_{2n+2}-2n-5\equiv 2F_{2n+2}^2-5\equiv 0\pmod{3},$$
we have
\begin{align}\label{eq1}
 3\mid \sum_{i=1}^nF_{2i}^2.
\end{align}
Again, since $\pi(2^{k+1})=3\cdot 2^{k}$ and $\pi(2^{k+1})\mid 2n$, by Theorem \ref{thm1.2}, we get $2^{k+1}\mid F_{2n}$ and $F_{2n+1}\equiv F_{2n+2}\equiv 1 \pmod{2^{k+1}}$. Since $2^{k+1}\mid 2n$, it follows that
$$3F_{2n+1}^2+2F_{2n+2}^2-6F_{2n}F_{2n+2}-2n-5\equiv 3+2-5\equiv 0 \pmod{2^{k+1}},$$
which implies that
\begin{align}\label{eq2}
    2^{k+1}\mid \sum_{i=1}^nF_{2i}^2.
\end{align}
Therefore, from (\ref{eq1}) and (\ref{eq2}), we get
$$\frac{1}{n} \sum_{i=1}^nF_{2i}^2\in \mathbb{N}.$$
Set $n=5\cdot11^k$ for $k\geq1$. Then since $\pi(11^k)\mid 10\cdot 11^{k-1}\mid 2n$, by Theorem \ref{thm1.2}, we get $11^{k}\mid F_{2n}$ and $F_{2n+1}\equiv F_{2n+2}\equiv 1 \pmod{11^{k}}$. Therefore, as $11^k\mid 2n$,
$$3F_{2n+1}^2+2F_{2n+2}^2-6F_{2n}F_{2n+2}-2n-5\equiv 3+2-5\equiv 0 \pmod{11^{k}},$$
and it follows that
\begin{align}\label{eq3}
    11^k\mid \sum_{i=1}^nF_{2i}^2.
\end{align}
Since $\pi(25)=100$ and $2n=10\cdot 11^k\equiv 10 \pmod{100}$, we get
$$F_{2n}\equiv F_{10}, \quad F_{2n+1}\equiv F_{11},\quad F_{2n+2}\equiv F_{12},\quad 2n\equiv 10\pmod{25}.$$
Using the above congruences, we have
$$3F_{2n+1}^2+2F_{2n+2}^2-6F_{2n}F_{2n+2}-2n-5\equiv 3F_{11}^2+2F_{12}^2-6F_{10}F_{12}-15=17700\equiv 0\pmod{25},$$
from which we obtain
\begin{align}\label{eq4}
 5\mid   \sum_{i=1}^nF_{2i}^2.
\end{align}
Therefore, from (\ref{eq3}) and (\ref{eq4}), we get $$\frac{1}{n}\sum_{i=1}^nF_{2i}^2\in \mathbb{N},$$
which proves (2). Similarly, take $n=3\cdot 2^{k+1}$ when $n$ is even and $n=5\cdot 11^k$ when $n$ is odd, with $$\sum_{i=1}^nF_{2i-1}^2=\frac{3F_{2n}^2+2F_{2n-1}^2-4F_{2n-2}F_{2n}+2n-2}{5},$$ we obtain the proof of $(3)$ by applying the same argument as in the preceding case.

From Proposition \ref{proposition2}, we have
   $$\sum_{i=1}^{n} F_{3i}^{2}=\frac{L_{6n+3}-4(-1)^n}{20}.$$
 Set $n=2^k$, where $k\geq1$. Since $\pi(2^{k+2})=3\cdot 2^{k+1}=6n$, we have $Q^{6n}\equiv I\pmod{2^{k+2}}$, and hence
$$L_{6n+m}\equiv L_m \pmod{2^{k+2}},\quad \text{for all integers } m.$$
In particular, since $L_3=4$ and $n$ is even, the numerator $L_{6n+3}-4(-1)^n=L_{6n+3}-4$ is divisible by $2^{k+2}$. It is also divisible by $20$, so it is divisible by $\mathrm{lcm}(20,2^{k+2})=5\cdot2^{k+2}$, and therefore
$$n=2^k\mid \sum_{i=1}^{n} F_{3i}^{2}.$$
Set $n=5^k$, where $k\geq 1$. Since $n$ is odd, $\displaystyle{\sum_{i=1}^{n} F_{3i}^{2}=(L_{6n+3}+4)/20}$, so it suffices to prove $L_{6n+3}\equiv -4\pmod{5^{k+1}}$. By Lemma \ref{lemma1.7}, we have $5^k\mid F_n$. By Proposition \ref{prop:basicid}(4) and (5), $L_n^2=5F_n^2-4$ and $2L_{n+1}=L_n+5F_n$, hence
$$2L_nL_{n+1}=L_n^2+5F_nL_n=5F_n^2-4+5F_nL_n\equiv -4\pmod{5^{k+1}},$$
so $L_nL_{n+1}\equiv -2\pmod{5^{k+1}}$. Using $L_{2n+1}=L_nL_{n+1}-(-1)^n$ (see \cite[p. 90]{r13}), we get $x:=L_{2n+1}\equiv -1\pmod{5^{k+1}}$. Finally, $L_{3m}=L_m^3-3(-1)^mL_m$ with $m=2n+1$ gives
$$L_{6n+3}=x^3+3x\equiv -1-3=-4\pmod{5^{k+1}}.$$
Hence $\displaystyle{5^k\mid \sum_{i=1}^{5^k}F_{3i}^2}$. This proves (4).
\end{proof}



\begin{proof}[Proof of Theorem \ref{theorem p}]

From Proposition \ref{proposition2}, we have

  $$\sum_{i=1}^{n} F_{3i-1}^{2}=\frac{L_{6n+1}+4(-1)^n-5}{20}, \quad \sum_{i=1}^{n} F_{3i-2}^{2}=\frac{L_{6n-1}+5-4(-1)^n}{20}.$$
Set $n=2^k$, where $k\geq1$. As in the proof of Theorem \ref{Theorem 5}(4), we have $L_{6n+m}\equiv L_m \pmod{2^{k+2}}$ for all integers $m$. Since $n$ is even, $L_1=1$ and $L_{-1}=-1$, we get
$$L_{6n+1}+4-5\equiv 0,\quad L_{6n-1}+5-4\equiv 0 \pmod{2^{k+2}}.$$
Both numerators are also divisible by $20$, so the two quotients are divisible by $2^k$. Consequently, we get
$$n\mid \sum_{i=1}^{n} F_{3i-1}^{2},\quad n\mid \sum_{i=1}^{n} F_{3i-2}^{2}.$$

 Let $p$ be an odd prime. For $p=5$ we have $\displaystyle{\sum_{i=1}^5F_{3i-1}^2=150517}$ and $\displaystyle{\sum_{i=1}^5F_{3i-2}^2=57493}$, neither of which is divisible by $5$. Let $p\neq 5$. Since $p$ is odd, it suffices to show that
$L_{6p+1}\not\equiv 9$ and $L_{6p-1}\not\equiv -9\pmod{p}$.

If $p\equiv \pm 2\pmod{5}$, then $\pi(p)\mid 2(p+1)\mid 6p+6$, so $L_{m+6p+6}\equiv L_m\pmod p$ for all $m$, and therefore
$$L_{6p+1}\equiv L_{-5}=-11 \not\equiv 9\pmod{p},\qquad L_{6p-1}\equiv L_{-7}=-29 \not\equiv -9\pmod{p},$$
because $p\nmid 20$.

If $p\equiv \pm 1\pmod{5}$, then $\pi(p)\mid p-1\mid 6p-6$, and therefore
$$L_{6p+1}\equiv L_{7}= 29\not\equiv 9 \pmod{p},\qquad L_{6p-1}\equiv L_5=11\not\equiv -9 \pmod{p}.$$
Therefore, for every odd prime $p$ we have
$$p\nmid \sum_{i=1}^{p} F_{3i-1}^{2},\quad p\nmid \sum_{i=1}^{p} F_{3i-2}^{2}.$$
This completes the proof.
\end{proof}

\begin{proof}[Proof of Theorem \ref{Theorem 4}]
By Proposition \ref{prop:basicid}(1),
we have
$$5\sum_{i=1}^n F_{r+id}^2=\sum_{i=1}^nL_{2r+2id}-2\sum_{i=1}^n (-1)^{r+id}=\sum_{i=1}^nL_{2r+2id}-2(-1)^r\sum_{i=1}^n (-1)^{id}.$$
Since $L_n = \tr (Q^n)$, we have
$$\sum_{i=1}^nL_{2r+2id}=\tr\biggl(\sum_{i=1}^nQ^{2r+2id}\biggr).$$
Now, let $n=8\cdot 3^k$. By Lemma \ref{lem:key} (applied with $2r$, $2d$ and $\lambda=1$), we have $\displaystyle{\sum_{i=1}^{n}Q^{2r+2id}\equiv 0\pmod n}$, and therefore
$$\sum_{i=1}^{8\cdot 3^{k}}L_{2r+2id}\equiv 0 \pmod{8\cdot 3^k}.$$ Further, note that
\[2\sum_{i=1}^{8\cdot 3^{k}} (-1)^{r+id}=2(-1)^r\sum_{i=1}^{8\cdot 3^{k}} (-1)^{id} = \begin{cases}
    0  &  \text{if } d \text{ is odd},\\
    2(-1)^r\, 8\cdot 3^{k} & \text{if } d \text{ is even.}
\end{cases}
\]
Consequently, $\displaystyle{5\sum_{i=1}^{8\cdot 3^k}F_{r+id}^2\equiv 0 \pmod{8\cdot 3^k}}$, and since $\gcd(5,8\cdot 3^k)=1$,
$$\sum_{i=1}^{8\cdot 3^k}F_{r+id}^2\equiv 0 \pmod{8\cdot 3^k}.$$

Now assume that $4\mid d$ and let $n=5^k$. Put $\displaystyle{T_q=\sum_{i=1}^{5}F_{r+id}^q}$ for $q=1,2,3$. First we show that $5\mid T_q$. By the proof of Lemma \ref{lem:lift5}, $Q^d\equiv I+cN\pmod 5$ for some integer $c$, where $N=\left(\begin{smallmatrix}4&3\\3&1\end{smallmatrix}\right)$, so $Q^{r+id}\equiv Q^r(I+icN)\pmod 5$, and hence $F_{r+id}\equiv a+ib\pmod 5$ for some integers $a,b$ independent of $i$. Since $\displaystyle{\sum_{i=1}^5i^j\equiv0\pmod5}$ for $j=1,2,3$, we get
$$T_q\equiv \sum_{i=1}^{5}(a+ib)^q\equiv 0\pmod 5\quad\text{where}\quad q=1,2,3.$$
Next, let $A=Q^{2d}$. By Lemma \ref{lem:lift5},
$$\sum_{i=1}^{5^k}L_{2r+2id}=\tr\Bigl(Q^{2r}\sum_{i=1}^{5^k}A^i\Bigr)\equiv 5^{k-1}\tr\Bigl(Q^{2r}\sum_{i=1}^{5}A^i\Bigr)=5^{k-1}\sum_{i=1}^{5}L_{2r+2id}\pmod{5^{k+1}}.$$
Since $d$ is even, $\displaystyle{\sum_{i=1}^{5}L_{2r+2id}=5T_2+10(-1)^r}$ and $\displaystyle{\sum_{i=1}^{5^k}(-1)^{id}=5^k}$. Therefore
$$5\sum_{i=1}^{5^k}F_{r+id}^2\equiv 5^{k-1}\bigl(5T_2+10(-1)^r\bigr)-2(-1)^r5^k=5^kT_2\equiv 0\pmod{5^{k+1}},$$
that is, $\displaystyle{5^k\mid \sum_{i=1}^{5^k}F_{r+id}^2}$. This completes the proof.
\end{proof}

\begin{proof}[Proof of Theorem \ref{Theorem 8}]
From Proposition \ref{Iden}, we have
$$\sum_{i=1}^nF_i^3=\frac{F_{n+2}^3-3F_{n+1}^3+3(-1)^nF_n+2}{4}.$$
Set $n=3\cdot 2^{k+2}$, where $k\geq 1$. Since $\pi(3)=8\mid n$, we have by Theorem \ref{thm1.2} that $3\mid F_n$ and $F_{n+1}\equiv 1 \pmod{3}$, and thus $F_{n+2}=F_{n+1}+F_n\equiv 1 \pmod{3}$. As $n$ is even, it follows that
$$F_{n+2}^3-3F_{n+1}^3+3F_n+2\equiv 1-0+0+2\equiv 0 \pmod{3}.$$
Therefore
\begin{align}\label{eq5}
3\mid \sum_{i=1}^nF_i^3.
\end{align}
Now $\pi(2^{k+3})=3\cdot 2^{k+2}=n$, using Theorem \ref{thm1.2}, we get $2^{k+3}\mid F_n$ and $F_{n+1}\equiv1\pmod{2^{k+3}}$, hence
$F_{n+2}\equiv 1 \pmod{2^{k+3}}$. Writing $F_{n+2}=1+2^{k+3}t$, we get $F_{n+2}^3\equiv1\pmod{2^{k+3}}$,
that is,
\begin{align}\label{eq6}
2(1-F_{n+2}^3)\equiv 0 \pmod{2^{k+4}}.
\end{align}
Since $12\mid n$, by Lemma \ref{lemma1.7}, we have $\nu_2(F_n)=\nu_2(n)+2=k+4$, hence $F_{n+1}=F_{n+2}-F_n\equiv F_{n+2}\pmod{2^{k+4}}$ and $3F_n\equiv0\pmod{2^{k+4}}$, which together with (\ref{eq6}) gives
$$F_{n+2}^3-3F_{n+1}^3+3F_n+2\equiv F_{n+2}^3-3F_{n+2}^3 +2= 2(1-F_{n+2}^3)\equiv 0 \pmod{2^{k+4}}.$$
Therefore
$$2^{k+2}\mid \sum_{i=1}^nF_i^3,$$
and this together with (\ref{eq5}) implies $n\mid \displaystyle{\sum_{i=1}^nF_i^3}$, since $\gcd(3,2^{k+2})=1$. This gives infinitely many even $n$ in (1).

For the odd integers, set $n=m-1$ in Lemma \ref{lem2.2}. Then $n$ is odd and, by Lemma \ref{Lemma luca}(1), we have $n\mid F_{n+1}$ and $F_n\equiv 1\pmod{n}$, and thus $F_{n+2}=F_{n+1}+F_n\equiv 1\pmod{n}$. Since $n$ is odd, $(-1)^n=-1$, and therefore
$$F_{n+2}^3-3F_{n+1}^3-3F_n+2\equiv 1-0-3+2\equiv 0 \pmod{n},$$
that is, $n\mid \displaystyle{\sum_{i=1}^nF_i^3}$, since $\gcd(n,4)=1$. As the integers $m-1$ of Lemma \ref{lem2.2} form an infinite set of odd integers (Remark \ref{remark1}), this proves (1).

Next, from \cite{newnew,r10}, we can write
$$\sum_{i=1}^{n} F_{2i}^{3}=\begin{cases}
\dfrac{1}{4}F_n^{2}L_{n+1}^{2}F_{n-1}L_{n+2}, & \text{if } n \text{ is even},\\[1ex]
\dfrac{1}{4}L_n^{2}F_{n+1}^{2}L_{n-1}F_{n+2}, & \text{if } n \text{ is odd}.
\end{cases}$$
Set $n=3\cdot 2^{k+2}$ for $k\geq1$. Since $\pi(3)\mid n$ and $\pi(2^{k+3})=n$, by Theorem \ref{thm1.2} we get $3\mid F_n$ and $2^{k+3}\mid F_n$, that is, $3\cdot 2^{k+3}\mid F_{n}$. Hence $3^2\cdot 2^{2k+6}\mid F_n^2$, and therefore
$$n=3\cdot 2^{k+2}\ \Big|\ \frac{1}{4}F_n^{2}L_{n+1}^{2}F_{n-1}L_{n+2}=\sum_{i=1}^{n} F_{2i}^{3}.$$
On the other hand, let $p$ be a prime with $p\equiv \pm 2 \pmod{5}$ (for instance $p\equiv 3\pmod{10}$). By Lemma \ref{lemma1.8}, $p\mid F_{p+1}$, and therefore
$$L_p^{2}F_{p+1}^{2}L_{p-1}F_{p+2}\equiv0 \pmod{p},$$
consequently $\displaystyle{p\mid \sum_{i=1}^pF_{2i}^3}$, as $p$ is odd. By Dirichlet's theorem there are infinitely many such primes $p$, and this proves (2).

From Proposition \ref{Iden}, we have
$$\sum_{i=1}^nF_{2i-1}^3=\frac{F_{2n}^3+3F_{2n}}{4}=\frac{F_{2n}(F_{2n}^2+3)}{4}.$$
We may take the even integers to be $n=3\cdot 2^{k+1}$, and the odd integers to be $n=5\cdot 11^k$, where $k\geq1$ is an integer. Indeed, if $n=3\cdot 2^{k+1}$, then $\pi(3)=8\mid 2n$ and $\pi(2^{k+3})=3\cdot 2^{k+2}=2n$, so $3\cdot 2^{k+3}\mid F_{2n}$ by Theorem \ref{thm1.2}, and hence $n=3\cdot2^{k+1}\mid F_{2n}/4$, which gives $n\mid \displaystyle{\sum_{i=1}^nF_{2i-1}^3}$. If $n=5\cdot 11^k$, then $5\mid F_{2n}$ by Lemma \ref{lemma1.7}, and $11^k\mid F_{2n}$ by Theorem \ref{thm1.2}, as $\pi(11^k)\mid 11^{k-1}\pi(11)=10\cdot11^{k-1}\mid 2n$, so $n\mid F_{2n}$, and since $n$ is odd, $n\mid \displaystyle{\sum_{i=1}^nF_{2i-1}^3}$. This proves (3).

Finally, let $p$ be an odd prime. For $p=5$, we have $\displaystyle{\sum_{i=1}^5F_i^3=162}$, which is not divisible by $5$. Let $p\neq5$. Since $p$ is odd, to prove $p\nmid \displaystyle{\sum_{i=1}^pF_i^3}$, it suffices to prove that $p\nmid F_{p+2}^3-3F_{p+1}^3-3F_p+2$. If $p\equiv \pm 1\pmod{5}$, we have from Lemma \ref{lemma1.8} that $p\mid F_{p-1}$ and $F_p\equiv 1 \pmod{p}$. Therefore $F_{p+1}=F_{p}+F_{p-1}\equiv 1 \pmod{p}$ and consequently $F_{p+2}=F_{p+1}+F_p\equiv 2\pmod{p}$, so that
 $$F_{p+2}^3-3F_{p+1}^3-3F_p+2\equiv 8-3-3+2\equiv  4\not\equiv0\pmod{p}.$$
If $p\equiv \pm 2\pmod{5}$, we have from Lemma \ref{lemma1.8} that $p\mid F_{p+1}$ and $F_p\equiv -1 \pmod{p}$. Therefore $F_{p+2}=F_{p+1}+F_p\equiv-1\pmod{p}$, and it follows that
 $$F_{p+2}^3-3F_{p+1}^3-3F_p+2\equiv -1-0+3+2\equiv 4\not\equiv0 \pmod{p}.$$
In both cases $4\not\equiv0\pmod p$ because $p$ is odd. Hence $p\nmid \displaystyle{\sum_{i=1}^pF_i^3}$ for every odd prime $p$, which completes the proof.
\end{proof}

\begin{proof}[Proof of Theorem \ref{Theorem 9}] From Proposition \ref{proposition3}, for even $n$ we have
$$\sum_{i=1}^nF_{3i}^3=\frac{1}{5}\cdot\frac{F_{9n+9}+F_{9n}-34}{76}-\frac{3}{5}\cdot\frac{F_{3n+1}-1}{2},$$
$$\sum_{i=1}^{n} F_{3i-1}^{3}=\frac{1}{5}\cdot\frac{10F_{9n+1}+2F_{9n}-10}{76}+\frac{3}{5}\cdot\frac{F_{3n}}{2}, $$
$$\sum_{i=1}^{n} F_{3i-2}^{3}=\frac{1}{5}\cdot\frac{14F_{9n}-6F_{9n+1}+6}{76}-\frac{3}{5}\cdot
\frac{F_{3n-1}-1}{2},$$
where all the fractions with denominators $76$ and $2$ are integers (the first ones are the integers $\displaystyle{\sum_{i\le n} F_{9i}}$, $\displaystyle{\sum_{i\le n} F_{9i-3}}$, $\displaystyle{\sum_{i\le n}F_{9i-6}}$, and the second ones are given by Proposition \ref{proposition1}).
Set $n=2^k$, where $k\geq2$ (the case $k=1$ is checked directly: for $n=2$ the three sums are $520$, $126$ and $28$, all even). Since $\gcd(n,5\cdot 19)=1$, it is enough to prove that the numerators over $76$ are divisible by $2^{k+2}$ and the numerators over $2$ are divisible by $2^{k+1}$.

Since $12\mid 3n$ and $12\mid 9n$, Lemma \ref{lemma1.7} gives $\nu_2(F_{3n})=\nu_2(F_{9n})=k+2$. Since $\pi(2^{k+1})=3\cdot 2^k=3n$, we have $Q^{3n}\equiv Q^{9n}\equiv I\pmod{2^{k+1}}$, hence
$$F_{3n+1}\equiv F_{3n-1}\equiv F_{9n+1}\equiv 1\pmod{2^{k+1}}.$$
This already shows $2^{k+1}\mid F_{3n+1}-1$, $2^{k+1}\mid F_{3n}$ and $2^{k+1}\mid F_{3n-1}-1$. Moreover,
\begin{align*}
10F_{9n+1}+2F_{9n}-10&=10(F_{9n+1}-1)+2F_{9n}\equiv 0\pmod{2^{k+2}},\\
14F_{9n}-6F_{9n+1}+6&=14F_{9n}-6(F_{9n+1}-1)\equiv 0\pmod{2^{k+2}}.
\end{align*}
Finally, $F_{3n+3}=2F_{3n+1}+F_{3n}\equiv 2\pmod{2^{k+2}}$, and since $3n+3$ is odd, Proposition \ref{prop:basicid}(2) gives
$$F_{9n+9}=5F_{3n+3}^3-3F_{3n+3}=F_{3n+3}(5F_{3n+3}^2-3)\equiv 2\cdot 17=34\pmod{2^{k+2}}.$$
Hence $2^{k+2}\mid F_{9n+9}+F_{9n}-34$. Consequently $n$ divides all three sums.
\end{proof}

\begin{proof}[Proof of Theorem \ref{Theorem 7}]
As $5F_n^3=F_{3n}-3(-1)^nF_n$, we have
\begin{equation}\label{eq:cube}
5\sum_{i=1}^n F_{r+id}^3=\sum_{i=1}^nF_{3r+3id}-3(-1)^r\sum_{i=1}^n \bigl((-1)^{d}\bigr)^{i}F_{r+id}.
\end{equation}
Let $n=8\cdot 3^{k}$. By Theorem \ref{thm:lambda} applied with $(3r,3d,\lambda=1)$ and with $(r,d,\lambda=(-1)^d)$, both sums on the right-hand side of \eqref{eq:cube} are divisible by $n$. Since $\gcd(5,n)=1$, we obtain
$$8\cdot 3^k\mid \sum_{i=1}^{8\cdot 3^k}F_{r+id}^3.$$

Now assume $4\mid d$ and let $n=5^k$, then $(-1)^d=1$. By Lemma \ref{lem:lift5} applied to $Q^{3d}$ and to $Q^d$ (extracting the $(1,2)$-entries after multiplying by $Q^{3r}$ and $Q^r$, respectively),
$$\sum_{i=1}^{5^k}F_{3r+3id}\equiv 5^{k-1}\sum_{i=1}^{5}F_{3r+3id},\qquad \sum_{i=1}^{5^k}F_{r+id}\equiv 5^{k-1}\sum_{i=1}^{5}F_{r+id}\pmod{5^{k+1}}.$$
Hence, by \eqref{eq:cube} (used for $n=5^k$ and for $n=5$),
$$5\sum_{i=1}^{5^k}F_{r+id}^3\equiv 5^{k-1}\cdot 5\sum_{i=1}^5F_{r+id}^3=5^{k}T_3\equiv 0\pmod{5^{k+1}},$$
where $T_3=\sum_{i=1}^{5}F_{r+id}^3$, and $T_3\equiv 0\pmod 5$ was shown in the proof of Theorem \ref{Theorem 4}. Thus $\displaystyle{5^k\mid \sum_{i=1}^{5^k}F_{r+id}^3}$.
\end{proof}

\begin{proof}[Proof of Theorem \ref{Theorem 4a}]
Since $\displaystyle{\sum_{i=1}^ni=\frac{n(n+1)}{2}}$ and $\displaystyle{\sum_{i=1}^nF_i^2=F_nF_{n+1}}$, it suffices to show that $n(n+1)\mid 2F_nF_{n+1}$.

Set $n=m$ in Lemma \ref{lem2.2}, then $n$ is even. By Remark \ref{remark1}(c), we have $n\mid F_{n/2}$, and hence $n\mid F_n$, as $F_{n/2}\mid F_n$. By Lemma \ref{Lemma luca}(2), we have $n+1\mid F_n$. Since $\gcd(n,n+1)=1$, we get $n(n+1)\mid F_n\mid 2F_nF_{n+1}$.

Now set $n=m-1$ in Lemma \ref{lem2.2}, then $n$ is odd. By Lemma \ref{Lemma luca}(1), we have $n\mid F_{n+1}$, and from Remark \ref{remark1}(c), we have $n+1=m\mid F_{m/2}\mid F_{m}=F_{n+1}$. Since $\gcd(n,n+1)=1$, we have $n(n+1)\mid F_{n+1}\mid 2F_nF_{n+1}$.

Since the integers $m$ (respectively $m-1$) of Lemma \ref{lem2.2} form an infinite set of even (respectively odd) positive integers (Remark \ref{remark1}), the theorem follows.
\end{proof}
\begin{proof}[Proof of Theorem \ref{Theorem 10}]
For $n=1,2$ we have $F_1^k\cdots F_n^k=1$, so the assertion holds. Let $n\geq3$ with $n\neq 6,12$, and let $p$ be a primitive prime divisor of $F_n$. Then $p$ divides none of $F_1,\ldots,F_{n-1}$, so $\nu_p(F_1^k\cdots F_n^k)=k\,\nu_p(F_n)$. We also use that $F_n<2^n$ for all $n\geq1$.

Let $k=1$. If $(F_1F_2 \cdots F_n)^{\frac{1}{n}} \in \mathbb{N}$, then $n\mid \nu_p(F_n)$, so we must have $p^{n}\mid F_n$ and therefore
$$2^{n}\leq p^{n}\leq F_n<2^{n},$$
but this is absurd.

Let $k\geq2$. If $(F_1^kF_2^k \cdots F_n^k)^{\frac{1}{n^k}} \in \mathbb{N}$, then $n^k\mid k\,\nu_p(F_n)=\nu_p(F_n^k)$, so we must have $p^{n^k}\mid F_n^k$ and therefore
$$2^{n^k}\leq p^{n^k}\leq F_n^k<2^{kn},$$
that is, $n^k<kn$, or equivalently $n^{k-1}<k$, which is absurd since $n^{k-1}\geq 2^{k-1}\geq k$.

Now for $n=6,12$, we have $F_1F_2\cdots F_6=2^4\cdot3\cdot5$ and $F_1F_2\cdots F_{12}=2^{9}\cdot3^{4}\cdot5^{2}\cdot 7\cdot 11\cdot13 \cdot 17\cdot 89$, and therefore
$$(F_1^kF_2^kF_3^kF_4^kF_5^kF_6^k)^{\frac{1}{6^k}}=(2^{4k}\cdot3^k\cdot5^k)^{\frac{1}{6^k}}\not\in \mathbb{N},$$
$$(F_1^kF_2^k\cdots F_{12}^k)^{\frac{1}{12^k}}=(2^{9k}\cdot3^{4k}\cdot5^{2k}\cdot 7^k\cdot 11^k\cdot13^k \cdot 17^k\cdot 89^k )^{\frac{1}{12^k}}\not\in \mathbb{N},$$
since the exponent of $3$ (respectively of $7$) is $k$, while $6^k\nmid k$ and $12^k\nmid k$ (as $6^k>k$ and $12^k>k$). This completes the proof.
\end{proof}

\noindent For integers $n,k\geq1$, we write
\begin{equation}\label{eq:Snk}
S_{n,k}=\sum_{i=1}^{n}\frac{1}{F_i^{k}}=\frac{U_{n,k}}{V_{n,k}},
\end{equation}
where $U_{n,k}$ and $V_{n,k}$ are coprime positive integers, that is, $U_{n,k}/V_{n,k}$ is $S_{n,k}$ in lowest terms. Since $\gcd(U_{n,k},V_{n,k})=1$, the number $n^k/S_{n,k}=n^kV_{n,k}/U_{n,k}$ (which for $k=1,2$ is the quantity in Theorem \ref{Theorem 11}) is an integer if and only if $U_{n,k}\mid n^k$. The proof of Theorem \ref{Theorem 11} is $2$-adic: for $n\geq6$ and $k\geq1$ we show that $\nu_2(V_{n,k})=k(h+3)$, where $h=\floor{\log_2(n/6)}$, and hence $U_{n,k}=S_{n,k}V_{n,k}\geq 2^{k(h+3)}S_{n,k}$. Since $n<12\cdot2^{h}$, this gives $U_{n,k}>n^k$, and hence $U_{n,k}\nmid n^k$, as soon as $S_{n,k}\geq(3/2)^k$. This holds for $k=1,2$, but fails for every $k\geq3$, because $S_{n,k}<2.3<(3/2)^3$ for all $n\geq1$ and $k\geq3$. So a different idea is needed for $k\geq3$ (see Conjecture \ref{conj:harm}).

\begin{proof}[Proof of Theorem \ref{Theorem 11}] Let $n\geq 6$ and $k\geq1$, and let $S_{n,k}$, $U_{n,k}$ and $V_{n,k}$ be as in \eqref{eq:Snk}. Let $h=\floor{\log_2 (n/6)}$, so that $6\cdot 2^h\leq n<12\cdot 2^h$, and set $R=6\cdot 2^{h}$. By Lemma \ref{lemma1.7}, $\nu_2(F_R)=h+3$. We claim that $\nu_2(F_r)<h+3$ for every $r\in\{1,\ldots,n\}$ with $r\neq R$. Indeed, if $3\nmid r$ then $\nu_2(F_r)=0$, if $r\equiv 3\pmod 6$ then $\nu_2(F_r)=1$, if $r\equiv 6\pmod{12}$ then $\nu_2(F_r)=3$, and $r\neq R$ forces $h\geq1$ in this case (if $h=0$, the only multiple of $6$ up to $n<12$ is $R=6$), finally, if $12\mid r$, then $\nu_2(F_r)=\nu_2(r)+2$, and $\nu_2(r)\geq h+1$ together with $3\mid r\leq n<12\cdot 2^h$ forces $r=3\cdot 2^{h+1}=R$. This proves the claim.

Consequently, $\nu_2(1/F_R^k)=-k(h+3)$ is strictly smaller than $\nu_2(1/F_r^k)$ for all $r\neq R$, and using the fact that for $x_1,x_2\in \mathbb{Q}$ with $\nu_2(x_1)\neq \nu_2(x_2)$ we have $\nu_2(x_1+x_2)=\min\{\nu_2(x_1),\nu_2(x_2)\}$, we get
$$\nu_2(S_{n,k})=\nu_2\Bigl(\frac{1}{F_{R}^k}\Bigr)=-k(h+3).$$
Since $\gcd(U_{n,k},V_{n,k})=1$, it follows that $\nu_2(V_{n,k})=k(h+3)$, so $V_{n,k}\geq 2^{k(h+3)}$. Now let $k\in\{1,2\}$. Then $S_{n,1}\geq S_{6,1}=\frac{379}{120}>3$ and $S_{n,2}\geq S_{6,2}>\frac94$. Therefore
$$U_{n,1}=S_{n,1}V_{n,1}>3\cdot 2^{h+3}=24\cdot 2^h>n,\qquad U_{n,2}=S_{n,2}V_{n,2}>\frac94\cdot 2^{2h+6}=144\cdot 4^h>n^2.$$
Now, if $\dfrac{n^k}{S_{n,k}}=\dfrac{n^kV_{n,k}}{U_{n,k}}\in \mathbb{N}$, then $U_{n,k}\mid n^k$ because $\gcd(U_{n,k},V_{n,k})=1$, in particular $U_{n,k}\leq n^k$, a contradiction. Hence $n^k/S_{n,k}\notin\mathbb{N}$ for all $n\geq 6$ and $k=1,2$.

It remains to check $n\leq 5$. For $k=1$, the values of $n/S_{n,1}$ for $n=1,\ldots,5$ are $1,1,\frac65,\frac{24}{17},\frac{150}{91}$, for $k=2$, the values of $n^2/S_{n,2}$ are $1,2,4,\frac{576}{85},\frac{22500}{2161}$. This proves the theorem.
\end{proof}

In the remaining proofs of this section we shall use repeatedly the criterion \eqref{eq:cong} for membership in $\cA$ and the congruences \eqref{eq:F1mod345}.

\begin{proof}[Proof of Theorem \ref{thm:Api}]
Let $n\geq2$ with $\pi(n)\mid n$. Since $\pi(n)$ is the order of $Q$ modulo $n$, we have $Q^n\equiv I\pmod n$, and therefore
$$Q^{n+2}=Q^{n}Q^{2}\equiv Q^2\pmod n.$$
Comparing the $(1,2)$-entries, we get $F_{n+2}\equiv F_2=1\pmod n$, that is, $n\in\cA$ by \eqref{eq:cong}.
Now let $n=2^a3^b5^c$ with $a\geq3$, $b\geq1$ and $c\geq0$. We have $\pi(2^a)=3\cdot2^{a-1}$, $\pi(3^b)\mid 3^{b-1}\pi(3)=8\cdot3^{b-1}$ (by property (1) of the Pisano period), and $\pi(5^c)=4\cdot5^c$ for $c\geq1$ (while $\pi(5^0)=\pi(1)=1$ divides $4$). Hence
$$\pi(n)=\mathrm{lcm}\bigl(\pi(2^a),\pi(3^b),\pi(5^c)\bigr)\ \Big|\ \mathrm{lcm}\bigl(3\cdot 2^{a-1},\,8\cdot 3^{b-1},\,4\cdot 5^c\bigr)=2^{\max\{a-1,3\}}\cdot3^{\max\{1,b-1\}}\cdot5^c.$$
Since $a\geq3$ and $b\geq1$, we have $\max\{a-1,3\}\leq a$ and $\max\{1,b-1\}\leq b$, thus $\pi(n)\mid n$, and $n\in\cA$ by what we have just shown. This proves (1).

Let $r\geq3$. The three integers
$$3\cdot2^r,\qquad 6\cdot2^r=3\cdot2^{r+1},\qquad 9\cdot2^r=3^2\cdot2^r$$
are of the form $2^a3^b$ with $(a,b)=(r,1)$, $(r+1,1)$ and $(r,2)$, respectively, where $a\geq3$ and $b\geq1$. Hence all three belong to $\cA$ by (1). Since $6\cdot2^r-3\cdot2^r=9\cdot2^r-6\cdot2^r=3\cdot2^r$, they form a $3$-term arithmetic progression. This proves (2).

Finally, let $k\geq3$, let $P\geq \max\{k,3\}$ be an integer and let $\mathcal{P}$ be the set of primes $p\leq P$. For $p\in\mathcal{P}$, every prime factor $q$ of $\pi(p)$ lies in $\mathcal{P}$: indeed $\pi(2)=3$, $\pi(3)=8$, $\pi(5)=20$, and for $p\geq7$, $\pi(p)$ divides $p-1$ or $2(p+1)$ (properties (2) and (3) of the Pisano period), whose prime factors are at most $p$. Let
$$c=\max\{\nu_q(\pi(p)) : p,q\in\mathcal{P}\},$$
and for an integer $E\geq \max\{c,\nu_2(k!)\}$ put
$$N=N_E=\prod_{p\in\mathcal{P}}p^{E}.$$
Then $k!\mid N$, because every prime factor $p$ of $k!$ lies in $\mathcal P$ and $\nu_p(k!)\leq\nu_2(k!)\leq E$. We claim that $\pi(jN)\mid jN$ for every $j\in\{1,\ldots,k\}$. Write
$$jN=\prod_{p\in\mathcal P}p^{E+f_p},\qquad f_p=\nu_p(j)$$
(all prime factors of $j$ lie in $\mathcal{P}$, since $j\leq k\leq P$). By property (1) of the Pisano period, $\pi(p^{E+f_p})\mid p^{E-1+f_p}\pi(p)$, and every prime factor of the right-hand side lies in $\mathcal P$. For a prime $q\in\mathcal P$ with $q\neq p$ we get
$$\nu_q\bigl(p^{E-1+f_p}\pi(p)\bigr)=\nu_q(\pi(p))\leq c\leq E\leq E+f_q,$$
while for $q=p$ we get
$$\nu_p\bigl(p^{E-1+f_p}\pi(p)\bigr)=E-1+f_p+\nu_p(\pi(p))\leq E+f_p,$$
because $\nu_p(\pi(p))=0$ for $p\neq5$ and $\nu_5(\pi(5))=\nu_5(20)=1$ (indeed $\pi(5^{E+f_5})=4\cdot5^{E+f_5}$). Hence $\pi(p^{E+f_p})\mid jN$ for every $p\in\mathcal P$, and therefore
$$\pi(jN)=\mathrm{lcm}_{p\in\mathcal P}\,\pi(p^{E+f_p})\ \Big|\ \prod_{p\in\mathcal P}p^{E+f_p}=jN.$$
By (1), $jN\in\cA$ for $j=1,\ldots,k$, that is, $N,2N,\ldots,kN$ is a $k$-term arithmetic progression contained in $\cA$. Since $E$ can be taken arbitrarily large and different values of $E$ give different values of $N$, we obtain infinitely many such progressions. This proves (3).
\end{proof}

\begin{proof}[Proof of Theorem \ref{thm:Aeven}]
Let $n\in\cA$ be even with $n>2$. By \eqref{eq:cong}, $n\mid F_{n+2}-1$, since $n$ is even, $F_{n+2}$ is odd. As $F_k$ is even if and only if $3\mid k$, we get $3\nmid n+2$, and therefore
$$n\equiv0\text{ or }2\pmod3.$$
Suppose first that $4\mid n$. Then \eqref{eq:cong} gives $F_{n+2}\equiv1\pmod4$, so by \eqref{eq:F1mod345}, $n+2\equiv1,2,5\pmod6$. Since $n+2$ is even, we must have $n+2\equiv2\pmod6$, and consequently $6\mid n$. Together with $4\mid n$, this gives $12\mid n$, that is,
$$n\equiv0\text{ or }12\pmod{24}.$$
If $n\equiv12\pmod{24}$, then $3\mid n$, and \eqref{eq:cong} gives $F_{n+2}\equiv1\pmod3$, however, $n+2\equiv14\equiv6\pmod8$, so $F_{n+2}\equiv F_6=8\equiv2\pmod3$, a contradiction. Hence $n\equiv0\pmod{24}$ in this case.

Suppose now that $4\nmid n$, so that $n\equiv2\pmod4$. If $n\equiv0\pmod3$, then $n\equiv6\pmod{12}$, i.e., $n\equiv6\text{ or }18\pmod{24}$, and $F_{n+2}\equiv1\pmod3$ by \eqref{eq:cong}. But if $n\equiv6\pmod{24}$, then $n+2\equiv0\pmod8$ and $F_{n+2}\equiv F_0=0\pmod3$, while if $n\equiv18\pmod{24}$, then $n+2\equiv4\pmod8$ and $F_{n+2}\equiv F_4=3\equiv0\pmod3$, both are impossible. Therefore $n\equiv2\pmod3$, and together with $n\equiv2\pmod4$, the Chinese Remainder Theorem gives $n\equiv2\pmod{12}$, that is, $n\equiv2\text{ or }14\pmod{24}$.

Combining the two cases, we conclude that $n\equiv0,2,\text{ or }14\pmod{24}$. This proves (1).

Now let $n\in\cA$ be odd. We first show that $5\nmid n$. Suppose that $5\mid n$. Then \eqref{eq:cong} implies $F_{n+2}\equiv1\pmod5$, so by \eqref{eq:F1mod345},
$$n+2\equiv1,2,8,19\pmod{20},\quad\text{i.e.,}\quad n\equiv19,0,6,17\pmod{20}.$$
Since $n$ is odd, $n\equiv17\text{ or }19\pmod{20}$.
On the other hand, since $5\mid n$ and $n$ is odd, $n\equiv5\text{ or }15\pmod{20}$, a contradiction. Hence $5\nmid n$.

Next we show that $n$ is not a prime power. Suppose, on the contrary, that $n=p^k$, where $p$ is an odd prime and $k\geq1$. By what we have just shown, $p\neq5$. By Lemma \ref{lem:Fpk}, we have $Q^{p^k}\equiv Q$ or $Q^{p^k}\equiv-Q^{-1}\pmod p$, and hence
$$Q^{n+2}=Q^{p^k}Q^2\equiv Q^3\quad\text{or}\quad Q^{n+2}\equiv-Q^{-1}Q^{2}=-Q\pmod p.$$
Comparing the $(1,2)$-entries, we get $F_{n+2}\equiv F_3=2$ or $F_{n+2}\equiv-F_1=-1\pmod p$. Since $p\geq3$, neither $2$ nor $-1$ is congruent to $1$ modulo $p$, whereas $n\in\cA$ gives $F_{n+2}\equiv1\pmod p$ by \eqref{eq:cong}. This contradiction shows that $n$ is not a prime power, in particular, no odd prime belongs to $\cA$.

Finally, suppose that $3\mid n$. Then \eqref{eq:cong} gives $F_{n+2}\equiv1\pmod3$, so by \eqref{eq:F1mod345}, $n+2\equiv1,2,7\pmod8$. Since $n+2$ is odd, $n+2\equiv1\text{ or }7\pmod8$, i.e.,
$$n\equiv7\text{ or }5\pmod8.$$
Since also $3\mid n$, the Chinese Remainder Theorem gives $n\equiv15$ or $21\pmod{24}$, in particular, $n\not\equiv3,9\pmod{24}$. This proves (2).
\end{proof}

\begin{remark}
Theorem \ref{thm:Aeven}(1) is best possible modulo $24$: we have $24\in\cA$, $758642=2\cdot17\cdot53\cdot421\in\cA$ with $758642\equiv2\pmod{24}$, and $10805846\in\cA$ with $10805846\equiv14\pmod{24}$. These are the smallest elements of $\cA$ in the classes $2$ and $14$ modulo $24$ (apart from $n=2$). In particular, all even elements $2<n<758642$ of $\cA$ are divisible by $24$ (see \cite{Codes} for Python codes).
\end{remark}

\begin{remark}
The smallest odd elements of $\cA$ are $77,319,323,1517,3021,4757,6479,\ldots$ (there are exactly $80$ odd elements below $2\cdot10^5$, see \cite{Codes} for Python codes)). All odd residue classes modulo $24$ other than $3$ and $9$ occur, for instance $1411699\in\cA$ and $1411699\equiv 19\pmod{24}$. More precisely, the smallest odd elements of $\cA$ in the classes $1,5,7,11,13,15,17,19,21,23$ modulo $24$ (see \cite{Codes} for Python codes) are
$$30049,\ 77,\ 319,\ 323,\ 179869,\ 82983,\ 78089,\ 1411699,\ 3021,\ 6479,$$
respectively, while the classes $3$ and $9$ are excluded by Theorem \ref{thm:Aeven}(2).
\end{remark}

\begin{proof}[Proof of Theorem \ref{thm:three}]
Since $A_F(3)=\frac{F_1+F_2+F_3}{3}=\frac{4}{3}\notin\mathbb{Z}$, we cannot have $n=1$, hence $n\geq2$, and so $n,n+1,n+2\in\cA$. Among $n$ and $n+1$ exactly one is even, call it $e$. If $e=n$, then $e+1=n+1\in\cA$, if $e=n+1$, then $e+1=n+2\in\cA$. In both cases $e,e+1\in\cA$. If $e=2$, then $n=2$ (since $n\geq2$), and $n+1=3\in\cA$, but $3\notin\cA$ by Theorem \ref{thm:Aeven}(2), as $3$ is an odd prime. Hence $e\geq4$.

By Theorem \ref{thm:Aeven}(1), $e\equiv0,2,14\pmod{24}$. If $e\equiv2\pmod{24}$, then $e+1\equiv3\pmod{24}$, so $e+1$ is an odd element of $\cA$ with $3\mid e+1$ and $e+1\equiv3\pmod 8$, contradicting Theorem \ref{thm:Aeven}(2). Therefore $e\equiv0$ or $14\pmod{24}$.

Finally, if $n$ is even, then $n=e$, and $n+2\in\cA$ is even with $n+2>2$ and $n+2\equiv2$ or $16\pmod{24}$, Theorem \ref{thm:Aeven}(1) excludes $16$, so $n\equiv0\pmod{24}$. If $n$ is odd, then $e=n+1$, so $n=e-1\equiv23$ or $13\pmod{24}$. Hence $n\equiv 0,13$ or $23\pmod{24}$, which completes the proof.
\end{proof}

\begin{proof}[Proof of Theorem \ref{thm:FmA}]
By \eqref{eq:cong}, $F_m\in\cA$ if and only if $F_{F_m+2}\equiv1\pmod{F_m}$.

Let $m\geq5$ be odd and suppose that $F_m\in\cA$. By Lemma \ref{lem:Fk1}(2), $F_{F_m+2}\equiv1\pmod{F_m}$ if and only if
$$F_m+2\equiv1,\ 2,\ 2m-2,\ -1\pmod{4m},\quad\text{i.e.,}\quad F_m\equiv-1,\ 0,\ 2m-4,\ -3\pmod{4m}.$$
Since $m$ is odd and $4\mid 4m$, these four residues are congruent to $3,0,2,1$ modulo $4$, respectively. The Fibonacci sequence modulo $4$ has period $6$ and $(F_k\bmod 4)_{k=0}^{5}=(0,1,1,2,3,1)$, for odd $m$ we have $m\equiv1,3,5\pmod 6$, and hence $F_m\equiv1\pmod4$ if $3\nmid m$, and $F_m\equiv2\pmod4$ if $3\mid m$. Therefore
$$F_m\equiv-3\pmod{4m}\ \text{ if } 3\nmid m,\qquad F_m\equiv2m-4\pmod{4m}\ \text{ if } 3\mid m.$$
Reducing modulo $m$, we get $m\mid F_m+3$ if $3\nmid m$, and $m\mid F_m+4$ if $3\mid m$. This proves (1).

We now prove (2). Suppose that $F_m\in\cA$. Since $m\geq5$ is odd, Theorem \ref{thm:FmA}(1) applies.

If $p\geq7$, then $3\nmid m$, so $m\mid F_m+3$, and in particular $p\mid F_{p^k}+3$. By Lemma \ref{lem:Fpk},
$$F_{p^k}+3\equiv\varepsilon^k+3\in\{2,4\}\pmod p,$$
which is not $0$ modulo $p$ since $p\geq7$. This is a contradiction.

If $p=5$, then $3\nmid m$, so $5\mid F_{5^k}+3$. But $5\mid F_{5^k}$ by Lemma \ref{lemma1.7}, so $F_{5^k}+3\equiv3\pmod5$, a contradiction.

If $p=3$, then, since $m=3^k\geq5$, we have $k\geq2$, so $9\mid m$, as $3\mid m$, Theorem \ref{thm:FmA}(1) gives $m\mid F_m+4$, and in particular $9\mid F_{3^k}+4$. Since $\pi(9)=24$, and $3^k\equiv 3\pmod{24}$ for odd $k\geq3$ and $3^k\equiv9\pmod{24}$ for even $k\geq2$, we get
$$F_{3^k}\equiv F_3=2\pmod 9\quad\text{or}\quad F_{3^k}\equiv F_9=34\equiv7\pmod9,$$
so $F_{3^k}+4\equiv6$ or $2\pmod 9$. Both are nonzero modulo $9$, a contradiction. This proves (2).

Let $m\geq4$ be even. By Lemma \ref{lem:Fk1}(1), $F_{F_m+2}\equiv1\pmod{F_m}$ if and only if
$$F_m+2\equiv1,\ 2,\ -1\pmod{2m},\quad\text{i.e.,}\quad F_m\equiv-1,\ 0,\ -3\pmod{2m}.$$
If $3\mid m$, then $F_m$ is even, since $2m$ is even, $F_m\equiv-1$ or $-3\pmod{2m}$ would force $F_m$ to be odd, so only the residue $0$ is possible. If $2m\mid F_m$, then $4\mid F_m$ (as $4\mid 2m$), so $6\mid m$ by Lemma \ref{lemma1.7}, then $3\mid 2m\mid F_m$, so $4\mid m$ by Lemma \ref{lemma1.7}, hence $12\mid m$, and clearly $m\mid F_m$. Conversely, if $12\mid m$ and $m\mid F_m$, then by Lemma \ref{lemma1.7},
$$\nu_2(F_m)=\nu_2(m)+2>\nu_2(m)+1=\nu_2(2m),$$
and $\nu_p(F_m)\geq\nu_p(m)=\nu_p(2m)$ for every odd prime $p$, thus $2m\mid F_m$. This proves (3).

If $12\mid m$ and $m\mid F_m$, then $F_m\equiv0\pmod{2m}$ by what we have just shown, and hence $F_m\in\cA$ by (3). Moreover, if $m=12\cdot 3^a5^b=2^2\cdot3^{a+1}\cdot5^b$, then by Lemma \ref{lemma1.7},
$$\nu_2(F_m)=\nu_2(m)+2,\qquad \nu_3(F_m)=\nu_3(m)+1,\qquad \nu_5(F_m)=\nu_5(m),$$
so $m\mid F_m$, and therefore $F_{12\cdot 3^a5^b}\in\cA$. This proves (4)(a).

Let $p\equiv17$ or $23\pmod{30}$. Then $p$ is an odd prime with $p\equiv5\pmod 6$, so $p\equiv2\pmod3$, and $p\equiv\pm2\pmod5$. By Theorem \ref{thm:FC}, $F_{F_m+2}\equiv1\pmod{F_m}$ for $m=2p$ and $m=4p$, that is, $F_{2p},F_{4p}\in\cA$ by \eqref{eq:cong}. (Equivalently, the proof of Theorem \ref{thm:FC} shows that $F_{2p}\equiv-1\pmod{4p}$ and $F_{4p}\equiv-3\pmod{8p}$, which are among the residues in (3).) By Dirichlet's theorem there are infinitely many such primes $p$. This proves (4)(b).
\end{proof}

In particular, $F_p\notin\cA$ for every prime $p\geq5$. For composite odd $m$, Theorem \ref{thm:FmA}(1) still gives $F_m\equiv-3$ (or $-4$ if $3\mid m$) modulo every prime $p\mid m$, but if $m$ has two distinct prime factors, the residue of $F_m$ modulo $p$ is no longer determined by Lemma \ref{lem:Fpk}, and the argument does not extend.

\begin{proof}[Proof of Theorem \ref{thm:density}]
Fix $r\geq0$ and $d\geq1$, put $\displaystyle{S(n)=\sum_{i=0}^{n-1}F_{r+id}}$, and let $y>\max\{5,F_d\}$. If $p>y$ is a prime, then $p\neq5$, and $\pi(p)\nmid d$, otherwise $p\mid F_{\pi(p)}\mid F_d$ (Theorem \ref{thm1.2}), which gives $p\leq F_d$. Hence Lemma \ref{lem:classes} applies to every prime $p>y$.

Let $z>y$. If $n\in\cB_{r,d}$ and $p\mid n$ for some prime $p\in(y,z]$, then $p\mid n\mid S(n)$, so $n\in E_p$. Consequently, for every $x\geq1$,
\begin{equation}\label{eq:densitysplit}
\#\{n\leq x:\ n\in\cB_{r,d}\}\leq\#\{n\leq x:\ p\nmid n\ \text{for all primes } p\in(y,z]\}+\sum_{y<p\leq z}\#\{n\leq x:\ n\in E_p\}.
\end{equation}
By Lemma \ref{lem:classes}, $E_p$ is contained in at most $4d$ residue classes modulo $p\,\pi(p)$, so
\begin{equation}\label{eq:Epcount}
\#\{n\leq x:\ n\in E_p\}\leq4d\Bigl(\frac{x}{p\,\pi(p)}+1\Bigr).
\end{equation}

Let $P=\prod_{y<p\le z}p$ and
$$\mathcal{S}=\{n\in\mathbb Z:\ p\nmid n \text{ for all primes } y<p\le z\}.$$
If $n\equiv m\pmod P$, then $n\equiv m\pmod p$ for every prime $p\in(y,z]$, so
$p\mid n$ if and only if $ p\mid m$. Hence $n\in \mathcal{S}$ if and only if $m\in \mathcal{S}$, i.e.\ $\mathcal{S}$ is periodic modulo $P$. It follows that $\mathcal{S}$ has density $\prod_{y<p\leq z}(1-1/p)$. Dividing \eqref{eq:densitysplit} by $x$, using \eqref{eq:Epcount} and letting $x\to\infty$, we obtain for the upper density $\overline{d}(\cB_{r,d})$ that
$$\overline{d}(\cB_{r,d})\leq\prod_{y<p\leq z}\Bigl(1-\frac1p\Bigr)+4d\sum_{p>y}\frac{1}{p\,\pi(p)}.$$
Letting $z\to\infty$, the product tends to $0$, because $\sum_p1/p$ diverges. Then, letting $y\to\infty$, the remaining sum tends to $0$ by Lemma \ref{lem:sumpi}. Hence $\overline{d}(\cB_{r,d})=0$, that is, $\cB_{r,d}$ has natural density zero.
For the quantitative bound, we use the upper bound sieve estimate (see \cite{r15})
$$\#\{n\leq x:\ p\nmid n\ \text{for all primes } p\in(y,z]\}\ll x\prod_{y<p\leq z}\Bigl(1-\frac1p\Bigr),$$
which holds uniformly for $2\leq y<z\leq x$. (The periodicity of $\mathcal{S}$ gives this bound only for $x\geq P$, and $P$ is much larger than $x$ for the choice of $z$ below.) Let $x$ be large, and take $z=\sqrt x$ and $y=\exp\sqrt{\log x}$ (so that $y>\max\{5,F_d\}$ once $x$ is large in terms of $d$). By Mertens' theorem (see \cite[Theorem 429]{Hardy}), $\prod_{y<p\leq z}(1-1/p)\ll \log y/\log z$, and therefore
$$\#\{n\leq x:\ p\nmid n\ \text{for all primes } p\in(y,z]\}\ll x\,\frac{\log y}{\log z}=\frac{2x}{\sqrt{\log x}}.$$
By \eqref{eq:Epcount} and Lemma \ref{lem:sumpi},
$$\sum_{y<p\leq z}\#\{n\leq x:\ n\in E_p\}\ll d\,\frac{x}{\log y}+d\,z=d\,\frac{x}{\sqrt{\log x}}+d\sqrt x.$$
Inserting both bounds into \eqref{eq:densitysplit}, we get $\#\{n\leq x:\ n\in\cB_{r,d}\}\ll_d x/\sqrt{\log x}$ for all large $x$, and hence for all $x\geq3$.
\end{proof}

\begin{proof}[Proof of Corollary \ref{cor:density}]
(1) Since $\displaystyle{\sum_{i=1}^nF_i=\sum_{i=0}^{n-1}F_{1+i}}$, we have $\cA\cup\{1\}=\cB_{1,1}$, and the claim follows from Theorem \ref{thm:density} with $r=d=1$.

(2) We have $\displaystyle{\sum_{i=1}^{n}F_{r+di}=\sum_{i=0}^{n-1}F_{(r+d)+di}}$ with $r+d\geq0$, so the set in question is $\cB_{r+d,d}$, which has density zero by Theorem \ref{thm:density}. The sums in Corollary \ref{Theorem 1} and Theorem \ref{Theorem 2} are of this form, for instance, $\displaystyle{\sum_{i=1}^nF_{2i-1}}$ and $\displaystyle{\sum_{i=1}^nF_{3i-2}}$ correspond to $(r,d)=(-1,2)$ and $(-2,3)$.
\end{proof}

\section{Data Availability}
The authors confirm that their manuscript has no associated data.
\section{Competing Interests}
The authors confirm that they have no competing interests.
\section{Use of AI-assisted Tools}
Claude was used to discuss a few mathematical problems, mainly used in reviewing the manuscript. The ideas, results, and most of the mathematical contributions are entirely the authors' own. All results were independently verified by hand, and the authors take full responsibility for the accuracy, originality, and integrity of the manuscript.

\end{document}